\input amstex
\input amsppt.sty

\input epsf \epsfverbosetrue 
 \magnification
980\vsize=21 true cm \hsize=16.5 true cm \voffset=1.1 true cm
\pageno=1 \NoRunningHeads \TagsOnRight

\def\p{\partial}
\def\ve{\varepsilon}
\def\f{\frac}
\def\na{\nabla}

\def\al{\alpha}
\def\t{\tilde}
\def\q{\quad}
\def\vp{\varphi}
\def\O{\Omega}
\def\th{\theta}
\def\g{\gamma}
\def\G{\Gamma}

\def\dl{\delta}

\def\o{\omega}

\def\ds{\displaystyle}
\topmatter
\topmatter \vskip 0.3 true cm \title{\bf On the existence and cusp
singularity of solutions to semilinear generalized Tricomi equations
with discontinuous initial data}
\endtitle
\endtopmatter
\document

\vskip 0.2 true cm \footnote""{* Ruan Zhuoping and Yin Huicheng were
supported by the NSFC (No.~10931007, No.~11025105), by the Priority
Academic Program Development of Jiangsu Higher Education
Institutions, and by the DFG via the Sino-German project ``Analysis
of PDEs and application.'' This research was carried out when Witt
Ingo visited Nanjing University in March of 2012,  and Yin Huicheng
was visiting the Mathematical Institute of the University of
G\"{o}ttingen in July-August of 2012.} \footnote""{** Ingo Witt was
partly supported by the DFG via the Sino-German project ``Analysis
of PDEs and application.'' } \vskip 0.3 true cm \centerline{Ruan,
Zhuoping$^{1,*}$; \qquad Witt, Ingo$^{2,**}$;\qquad Yin,
Huicheng$^{1,*}$} \vskip 0.5 true cm {1. Department of Mathematics
and IMS, Nanjing University, Nanjing 210093, P.R.~China.}\vskip 0.2
true cm {2. Mathematical Institute, University of G\"{o}ttingen,
Bunsenstr.~3-5, D-37073 G\"{o}ttingen, Germany.} \vskip 0.3 true
cm

\centerline {\bf Abstract} \vskip 0.3 true cm In this paper, we are
concerned with the local existence and singularity structure of low
regularity solutions to the semilinear generalized Tricomi equation
$\p_t^2u-t^m\Delta u=f(t,x,u)$ with typical discontinuous initial data
$(u(0,x), \p_tu(0,x))=(0, \vp(x))$; here $m\in\Bbb N$, $x=(x_1, ...,
x_n)$, $n\ge 2$, and $f(t,x,u)$ is $C^{\infty}$ smooth in its
arguments. When the initial data $\vp(x)$ is a homogeneous function of
degree zero or a piecewise smooth function singular along the
hyperplane $\{t=x_1=0\}$, it is shown that the local solution
$u(t,x)\in L^{\infty}([0,T]\times\Bbb R^n)$ exists and is $C^{\infty}$
away from the forward cuspidal cone $\Gamma_0=\bigl\{(t,x)\colon t>0, \,
|x|^2=\ds\f{4t^{m+2}}{(m+2)^2}\bigr\}$ and the characteristic cuspidal
wedge $\G_1^{\pm}=\bigl\{(t,x)\colon t>0, x_1=\pm
\ds\f{2t^{\f{m}{2}+1}}{m+2}\bigr\}$, respectively. On the other hand, 
for $n=2$ and piecewise smooth initial data $\vp(x)$ singular along
the two straight lines $\{t=x_1=0\}$ and $\{t=x_2=0\}$, we establish
the local existence of a solution $u(t,x)\in L^{\infty}([0,
T]\times\Bbb R^2)\cap C([0, T], H^{\f{m+6}{2(m+2)}-}(\Bbb R^2))$ and
show further that $u(t,x)\not\in C^2((0,T]\times\Bbb
R^2\setminus(\G_0\cup\G_1^{\pm}\cup\G_2^{\pm}))$ in general due to the
degenerate character of the equation under study; here
$\G_2^{\pm}=\bigl\{(t,x)\colon t>0,
x_2=\pm\ds\f{2t^{\f{m}{2}+1}}{m+2}\bigr\}$. This is an essential
difference to the well-known result for solutions $v(t,x)\in
C^{\infty}(\Bbb R^+\times\Bbb R^2\setminus
(\Sigma_0\cup\Sigma_1^{\pm}\cup \Sigma_2^{\pm}))$ to the 2-D
semilinear wave equation $\p_t^2v-\Delta v=f(t,x,v)$ with
$(v(0,x), \p_tv(0,x))=(0, \vp(x))$, where $\Sigma_0=\{t=|x|\}$,
$\Sigma_1^{\pm}=\{t=\pm x_1\}$, and $\Sigma_2^{\pm}=\{t=\pm x_2\}$.

\vskip 0.3 true cm

\vskip 0.3 true cm

\vskip 0.2 true cm

\vskip 0.4 true cm \centerline{\bf $\S 1.$ Introduction} \vskip
0.2 true cm

In this paper, we will study the local existence and the singularity
structure of  low regularity solution to the following
$n$-dimensional semilinear generalized Tricomi equation
$$
\cases &\p_t^2u-t^m\Delta u=f(t,x,u),\qquad (t,x)\in [0,
+\infty)\times\Bbb R^n,\\
&u(0,x)=0,\quad \p_tu(0,x)=\vp(x),
\endcases\tag1.1
$$
where $m\in\Bbb N$, $x=(x_1, ..., x_n)$, $n\ge 2$,
$\Delta=\ds\sum_{i=1}^n\p_i^2$, $f(t,x,u)$ is $C^{\infty}$ smooth on
its arguments and has a compact support on the variable $x$, and the
typical discontinuous initial data $\vp(x)$ satisfies one of the
assumptions:

$(A_1)$ $\vp(x)=g(x,\ds\f{x}{|x|})$, here $g(x,y)\in
C^{\infty}(\Bbb R^n\times\Bbb R^n)$ and has a compact support in
$B(0,1)\times B(0, 2)$;

$(A_2)$ $\vp(x)=\cases \vp_1(x)\qquad\text{for
$x_1<0$},\\
\vp_2(x)\qquad\text{for
$x_1>0$},\\
\endcases$ with $\vp_1(x), \vp_2(x)\in C_0^{\infty}(\Bbb R^n)$ and $\vp_1(0)\not=\vp_2(0)$;

$(A_3)$ For $n=2$, $\vp(x)=\cases \psi_1(x)\qquad\text{for
$x_1>0$, $x_2>0$},\\
\psi_2(x)\qquad\text{for
$x_1<0$, $x_2>0$},\\
\psi_3(x)\qquad\text{for
$x_1<0$, $x_2<0$},\\
\psi_4(x)\qquad\text{for
$x_1>0$, $x_2<0$},\\
\endcases$ with $\psi_i(x)\in C_0^{\infty}(\Bbb R^n) (1\le i\le 4)$
and $\psi_i(0)\not=\psi_j(0)$ for some $i\not=j$ ($1\le i< j\le 4$).

It is noted that $\vp(x)=\psi(x)\ds\f{x_1}{|x|}$ with $\psi(x)\in C_0^{\infty}(B(0, 1))$
is a special function satisfying $(A_1)$, which has a singularity at the origin.

Under the assumptions $(A_1)-(A_3)$, we now state the main results in this
paper.

\vskip .2cm

 {\bf Theorem 1.1.} {\it There exists a constant $T>0$
such that

(i) Under the condition $(A_1)$, (1.1) has a unique solution
$u(t,x)\in C([0, T],$ $H^{\f{n}{2}+\f{2}{m+2}-}(\Bbb R^n))\cap C((0,
T],$
\newline
$H^{\f{n}{2}+\f{m+4}{2(m+2)}-}(\Bbb R^n))\cap C^1([0, T],
H^{\f{n}{2}-\f{m}{2(m+2)}-}(\Bbb R^n))$ and $u(t,x)\in
C^{\infty}((0, T]\times\Bbb R^n\setminus\G_0)$, here $\G_0=\{(t,x):
t>0, |x|^2=\ds\f{4t^{m+2}}{(m+2)^2}\}$.

(ii) Under the condition $(A_2)$, (1.1) has a unique solution
$u(t,x)\in L^{\infty}([0, T]\times\Bbb R^n)\cap C([0, T],$
$H^{\f{m+6}{2(m+2)}-}(\Bbb R^n))$ $\cap C((0, T],
H^{\f{m+3}{m+2}-}(\Bbb R^n))\cap C^1([0, T], H^{\f{1}{m+2}-}(\Bbb
R^n))$ and $u(t,x)\in C^{\infty}((0, T]\times\Bbb
R^n\setminus\G_1^+\cup\G_1^-)$, here $\G_1^{\pm}=\{(t,x): t>0,
x_1=\pm\ds\f{2t^{\f{m}{2}+1}}{m+2}\}$.

(iii) For $n=2$, under the condition $(A_3)$,   if $m\le 9$, then
(1.1) has a unique solution $u(t,x)\in L^{\infty}([0, T]\times\Bbb
R^2)\cap C([0, T],$ $H^{\f{m+6}{2(m+2)}-}(\Bbb R^2))\cap C((0, T],
H^{\f{m+3}{m+2}-}(\Bbb R^2))\cap C^1([0, T], H^{\f{1}{m+2}-}(\Bbb
R^2))$. Moreover, in the general case, $u(t,x)\not\in
C^2((0,T]\times\Bbb R^2\setminus\G_0\cup\G_1^{\pm}\cup\G_2^{\pm})$,
here $\G_0$ and $\G_1^{\pm}$ have been defined in (i) and (ii)
respectively, and $\G_2^{\pm}=\{(t,x): t>0, x_2=\pm
\ds\f{2t^{\f{m}{2}+1}}{m+2}\}$.}

\vskip .2cm

{\bf Remark 1.1.} {\it In order to prove the $C^{\infty}$ property
of solution in  Theorem 1.1.(i) and (ii), we will show that the
solution of (1.1) is conormal with respect to the cusp
characteristic conic surface $\G_0$ or the cusp characteristic
surfaces $\G_1^{\pm}$ respectively in $\S 6$ below. And the
definitions of conormal spaces will be given in $\S 4$.}

\vskip .2cm

{\bf Remark 1.2.} {\it Since we only focus on  the local existence
of solution in Theorem 1.1,  it does not lose the generality  that
the initial data $\varphi(x)$ in $(A_1)-(A_3)$ are assumed to be
 compactly supported. In addition, the initial data $(u(0,x),
\p_tu(0,x))=(0, \vp(x))$ in (1.1) can be   replaced by the general
forms $(u(0,x), \p_tu(0,x))=(\phi(x), \vp(x))$, where $D_x\phi(x)$ satisfies $(A_1)$ when $\vp(x)$ satisfies $(A_1)$,
$\phi(x)$ is $C^1$ piecewise smooth along $\{t=x_1=0\}$ when
$\vp(x)$ satisfies $(A_2)$, and $\phi(x)$ is $C^1$ piecewise smooth
along the lines $\{t=x_1=0\}$ and $\{t=x_2=0\}$  when $\vp(x)$
satisfies $(A_3)$, respectively.}

\vskip .2cm

{\bf Remark 1.3.} {\it The initial data problem (1.1) under the
assumptions $(A_2)$ and $(A_3)$ is actually a special case of the
multidimensional generalized Riemann problem for the second order
semilinear degenerate hyperbolic equations. For the semilinear
$N\times N$ {\bf strictly hyperbolic} systems of the form
$\p_tU+\ds\sum_{j=1}^nA_j(t,x)\p_jU=F(t,x,U)$ with the piecewise
smooth or conormal initial data along some hypersurface
$\Delta_0\subset\{(t,x): t=0, x\in\Bbb R^n\}$ (including the Riemann
discontinuous initial data), the authors in [19-20] have established
the local well-posedness of piecewise smooth or bounded conormal
solution with respect to the $N$ pairwise transverse characteristic
surfaces $\Sigma_j$ passing through $\Delta_0$. With respect to the
Riemann problem of higher order semilinear degenerate hyperbolic
equations, we will establish the related results in our forthcoming
paper.}

\vskip .2cm

{\bf Remark 1.4.} {\it The reason that we pose the restriction on
$m\le 9$ in (iii) of Theorem 1.1 is due to the requirement for
utilizing the Sobolev's imbedding theorem to derive the boundedness
of solution (one can see details in (5.7) of $\S 5$ below),
otherwise, it seems that we have to add some other conditions on the
nonlinear function $f(t,x,u)$ since the solution $w(t,x)\in
L^{\infty}([0, T]\times\Bbb R^2)$ does not hold even if $w(t,x)$
satisfies a linear equation $\p_t^2w-t^m\Delta w=g(t,x)$ with
$(w(0,x), \p_tw(0,x))=(0, 0)$ and $g(t, x)\in L^{\infty}([0,
T]\times\Bbb R^2)\cap L^{\infty}([0, T], H^s(\Bbb R^2))$ with $0\le
s<1$ for large $m$. Firstly, this can be roughly seen from the
following explicit formula of $w(t,x)$ in Theorem 3.4 of [24]:
$$
\align
w(t,x)&=\f{1}{\pi}(\f{4}{m+2})^{\f{m}{m+2}}\int_0^td\tau\int_0^{\phi(t)-\phi(\tau)}dr_1\biggl(\p_r
\int_{B(x,r)}\f{g(\tau,
y)}{\sqrt{r^2-|x-y|^2}}dy\biggr)\biggr|_{r=r_1}
(r_1+\phi(t)+\phi(\tau))^{-\g}\\
&\times(\phi(\tau)-r_1+\phi(t))^{-\g}F\biggl(\g,\g;1;
\f{(-r_1+\phi(t)-\phi(\tau))(-r_1-\phi(t)+\phi(\tau))}{(-r_1+\phi(t)+\phi(\tau))(-r_1-\phi(t)-\phi(\tau))}\biggr),
\endalign$$
where $\phi(t)=\ds\f{2t^{\f{m}{2}+1}}{m+2}$, $\g=\ds\f{m}{2(m+2)}$,
and $F(a,b;c;z)$ is the hypergeometric function. It is noted that
$\p_r\bigl(\int_{B(x,r)}\ds\f{g(\tau,
y)}{\sqrt{r^2-|x-y|^2}}dy\bigr)=\p_r\bigl(\int_0^r\int_0^{2\pi}g(\tau,
x_1+\sqrt{r^2-q^2}cos\th, x_2+\sqrt{r^2-q^2}sin\th)dqd\th\bigr)$
holds  and thus the $L^{\infty}$ property of $w(t,x)$ is closely
related to the integrability of the first order derivatives of
$g(t,x)$, which is different from the case in 2-D linear wave
equation. On the other hand, the regularity of $w(t,x)$ is in $
C([0, T], H^{s+\f{2}{m+2}}(\Bbb R^2))\not\subset L^{\infty}([0,
T]\times\Bbb R^2)$ for large $m$ by Proposition 3.3 below and
Sobolev's imbedding theorem.}

\vskip .2cm

{\bf Remark 1.5.} {\it By $u(t,x)\not\in C^2((0,T]\times\Bbb
R^2\setminus\G_0\cup\G_1^{\pm}\cup\G_2^{\pm})$ in Theorem 1.1.(iii),
we know that there exists an essential difference on the regularity
of solutions between the degenerate hyperbolic equation and strictly
hyperbolic equation with the same initial data in $(A_3)$ since
$v(t,x)\in C^{\infty}(\Bbb R^+\times\Bbb
R^2\setminus\Sigma_0\cup\Sigma_1^{\pm}\cup\Sigma_2^{\pm})$ will hold
true if  $v(t,x)$ is a  solution to the 2-D semilinear wave equation
$\p_t^2v-\Delta v=f(t,x,v)$ with $(v(0,x), \p_tv(0,x))=(0, \vp(x))$,
where $\Sigma_0=\{t=|x|\}$, $\Sigma_1^{\pm}=\{t=\pm x_1\}$, and
$\Sigma_2^{\pm}=\{t=\pm x_2\}$. The latter well-known result was
established in the references [1-2], [7-9] and [18] respectively
under some various assumptions.}

For $m=1, n=1$ and $f(t,x,u)\equiv 0$, the equation
in (1.1) becomes the classical
Tricomi equation which arises in transonic gas dynamics and has been
extensively investigated in bounded domain with suitable boundary
conditions from various viewpoints (one can see [4], [17], [21-22]
and the references therein). For $m=1$ and $n=2$, with respect to
the equation $\p_t^2u-t\triangle u=f(t,x,u)$ together with the
initial data of higher $H^s(\Bbb R^n)-$regularity $(s>\f{n}{2})$, M.Beals in [3] show the local
existence of solution $u\in C([0, T], H^s(\Bbb R^n))\cap C^1([0, T],
H^{s-\f{5}{6}}(\Bbb R^n))\cap C^2([0, T], H^{s-\f{11}{6}}(\Bbb
R^n))$ for some $T>0$ under the  crucial assumption that the
support of $f(t,x,u)$ on the variable $t$ lies in $\{t\ge 0\}$.
Meanwhile, the conormal regularity of $H^s(\Bbb R^n)$ solution
$u(t,x)$ with respect to the characteristic surfaces $x_1=\pm
\ds\f{2}{3}t^{\f{3}{2}}$ is also established in [3]. With respect to
more general nonlinear degenerate hyperbolic equations with higher
order regularities, the authors in [10-11] studied the local
existence and the propagation of weak singularity of {\bf classical
solution}. For the linear degenerate hyperbolic equations with
suitable initial data, so far there have existed some interesting
results on the regularities of solution when Levi's conditions are
posed (one can see [12], [14-15] and the references therein). In the
present paper, we focus  on the low regularity solution problem for
the second order semilinear degenerate equation with no much
restrictions on the nonlinear function $f(t,x,u)$ in (1.1) and typical discontinuous
initial data.

We now comment on the proof of Theorem 1.1. In order to prove the
local existence of solution to (1.1) with the low regularity, at
first we should establish the local $L^{\infty}$  property of
solution $v(t,x)$ to the linear problem $\p_t^2v-t^m\Delta v=F(t,x)$
with $(v(0,x),\p_tv(0,x))=(\vp_0(x), \vp_1(x))$ so that the
composition function $f(t,x,v)$ makes sense. In this process, we
have to make full use of the special structure of the piecewise
smooth initial data and the explicit expression of solution $v(t,x)$
established in [23-24] since we can not  apply for the Sobolev
imbedding theorem directly to obtain $v(t,x)\in L_{loc}^{\infty}$
due to its low regularity (for examples, in the cases of
$(A_2)-(A_3)$, the initial data are only in $H^{\f12-}(\Bbb R^n)$).
Based on such $L^{\infty}-$estimates, together with the Fourier
analysis method and the theory of confluent hypergeometric
functions, we can construct a suitable nonlinear mapping related to
the problem (1.1) and further  show that such a mapping admits a
fixed point in the space $L^{\infty}([0, T]\times\Bbb R^n)\cap C([0,
T], H^{s_0}(\Bbb R^n)$ for suitable $T>0$ and some number $s_0>0$,
and then the local solvability of (1.1) can be shown. Next, we are
concerned with the singularity structures of solution $u(t,x)$ of
(1.1). It is noted that the initial data are suitably conormal under
the assumptions $(A_1)$ and $(A_2)$, namely, $\ds\Pi_{1\le i,j\le
n}(x_i\p_j)^{k_{ij}}\vp(x)\in H^{\f{n}{2}-}(\Bbb R^n)$ for any
$k_{ij}\in\Bbb N\cup\{0\}$ in the case of $(A_1)$, and
$(x_1\p_1)^{k_1}\ds\Pi_{2\le i\le n}\p_i^{k_i}\vp(x)\in
H^{\f{1}{2}-}(\Bbb R^n)$ for any $k_i\in\Bbb N\cup\{0\}$ ($i=1, ...,
n$) in the case of $(A_2)$, then we intend to use the commutator
arguments as in [5-6] to prove  the conormality of solution $u(t,x)$
to (1.1). However, due to the cusp singularities of surfaces $\G_0,
\G_1^{\pm}$ together with the degeneracy of equation, it seems that
it is difficult to choose the smooth vector fields $\{Z_1, ...,
Z_k\}$ tangent to $\G_0$ or $\G_1^{\pm}$ as in [5-6] to define the
conormal space and take the related analysis on the commutators
$[\p_t^2-t^m\Delta, Z_1^{l_1}\cdot\cdot\cdot Z_k^{l_k}]$ since this
will lead to the violation of Levi's condition and bring the loss of
regularity of $Z_1^{l_1}\cdot\cdot\cdot Z_k^{l_k}u$ (more detailed
explanations can be found in $\S 4$ below). To overcome this
difficulty, motivated by [2-3] and [18], we will choose the
nonsmooth vector fields and try to find the extra regularity relations provided by
the operator itself and some parts of vector fields to yield full
conormal regularity of $u(t,x)$ together with the regularity theory
of second order elliptic equation and further complete the proof of
Theorem 1.1.(i) and (ii), here we point out that it is nontrivial
to find such crucial regularity relations. On the other hand, in the case of $n=2$
and assumption $(A_3)$, due to the lack of the strong Huyghen's
principle, we can derive that the solution $u(t,x)\not\in
C^2((0,T]\times\Bbb R^2\setminus\G_0\cup\G_1^{\pm}\cup\G_2^{\pm})$
of (1.1), which yields a  different phenomenon from that
in the case of second order strict hyperbolic equation as pointed
out in Remark 1.5.

This paper is organized as follows. In $\S 2$, for later uses,
we will give some preliminary results on the regularities of initial
data $\vp(x)$ in various assumptions $(A_1)-(A_3)$ and establish the
$L^{\infty}$ property of solution to the related linear problem. In
$\S 3$, by the partial Fourier-transformation, we can change the
linear generalized Tricomi equation into a confluent hypergeometric
equation, and then some weighted Sobolev regularity estimates near
$\{t=0\}$ are derived. In $\S 4$, the required conormal spaces are
defined and some crucial commutator relations are given. In $\S 5$,
based on the results in $\S 2$-$\S 3$, the local solvability of
(1.1) is established. In $\S 6$, we complete the proof on Theorem
1.1 by utilizing the concepts of conormal spaces and commutator
relations  in $\S 4$ and taking some analogous analysis in Lemma 2.4
of $\S 2$ respectively.

In this paper, we will use the following notation:
$$H^{s-}(\Bbb R^n)=\{w(x): \ w(x)\in H^{s-\dl}(\Bbb R^n) \quad \text {for any fixed constant } \dl>0. \}$$

\vskip 0.3 true cm \centerline{\bf $\S 2$. Some preliminaries}
\vskip 0.3 true cm

In this section, we will give some basic lemmas on the regularities
of initial data $\vp(x)$ in the assumptions $(A_1)-(A_3)$ and
establish some $L^{\infty}$ property of solution to the linear
problem $\p_t^2u-t^m\Delta u=f(t,x)$ with suitably piecewise smooth initial data.

With respect to the functions $\vp(x)$ given in $(A_1)-(A_3)$ of $\S
1$, we have the following regularities in Sobolev space.

{\bf Lemma 2.1.} {\it (i) If $\vp(x)=g(x,\ds\f{x}{|x|})$, here
$g(x,y)\in C^{\infty}(\Bbb R^n\times\Bbb R^n)$ and has a compact
support in $B(0, 1)\times B(0, 2)$, then $\vp(x)\in
H^{\f{n}{2}-}(\Bbb R^n)$.

(ii) If $n=2$ and $\vp(x)=\cases \psi_1(x)\qquad\text{for
$x_1>0$, $x_2>0$},\\
\psi_2(x)\qquad\text{for
$x_1<0$, $x_2>0$},\\
\psi_3(x)\qquad\text{for
$x_1<0$, $x_2<0$},\\
\psi_4(x)\qquad\text{for
$x_1>0$, $x_2<0$},\\
\endcases$ where $\psi_i(x)\in C_0^{\infty}(\Bbb R^2) (1\le i\le 4)$, then

$\vp(x)\in H^{\f12-}(\Bbb R^2)$.

(iii) If $\vp(x)=\cases \vp_1(x)\qquad\text{for
$x_1<0$},\\
\vp_2(x)\qquad\text{for
$x_1>0$},\\
\endcases$,  where $\vp_1(x), \vp_2(x)\in C_0^{\infty}(\Bbb R^n)$, then
$\vp(x)\in H^{\f 1 2-}(\Bbb R^n)$ and $x_1\vp(x)\in H^{\f 3 2-}(\Bbb
R^n)$.}

\vskip .2cm

 {\bf Proof.} (i)  It follows from a direct computation that
$$|\p_x^{\al}\vp(x)|\le C_{\al}|x|^{-|\al|}.\tag2.1$$

Since $\vp$ is integrable on $\Bbb R^n$, we have that $\hat{\vp}$,
the Fourier transform of $\vp$, is continuous on $\Bbb R^n$, which
implies that  $(1+|\xi|)^{\f{n}{2}-\dl}\hat{\vp}(\xi)\in
L^2(\{|\xi|\le 1\})$.

For $|\xi|>1$, we decompose $\hat{\vp}$ into two parts
$$
\align
\hat{\vp}(\xi)&=\int_{|x|<\f{1}{|\xi|}}e^{-ix \cdot\xi}\vp(x)dx+\int_{\f{1}{|\xi| } \le |x|\le 1}e^{-ix\cdot\xi}\vp(x)dx\\
&=I+II=I+\ds\sum_{\ell=1}^n\chi_\ell(\xi)II,\tag2.2
\endalign
$$
where $\{\chi_\ell\}_{\ell=1}^n$ is a $C^{\infty}$ conic
decomposition of unity corresponding to the domain $\{\xi\in\Bbb
R^n: |\xi|\ge 1\}$, moreover $\xi_\ell\not=0$ in $supp\chi_\ell$.

Obviously, the term $I$ can be dominated by the multiplier of
$|\xi|^{-n}$.

On the other hand, for any $1 \le \ell \le n,$
$$
\align
\chi_\ell(\xi) II&=\ds\f{\chi_\ell(\xi)}{|\xi|^n}\int_{1\le |x|\le |\xi|}
e^{-ix\cdot\f{\xi}{|\xi|}}g(\f{x}{|\xi|},\f{x}{|x|})dx\\
&=\ds\f{\chi_\ell(\xi)}{|\xi|^n}\ds\f{1}{i\f{\xi_\ell}{|\xi|}}\int_{1\le |x|\le |\xi|}e^{-ix\cdot\f{\xi}{|\xi|}}\p_\ell\biggl(g(\f{x}{|\xi|},\f{x}{|x|})\biggr)dx\\
&\quad +
\ds\f{\chi_\ell(\xi)}{|\xi|^n}\ds\f{1}{-i\f{\xi_\ell}{|\xi|}}\int_{|x|=1}e^{-ix
\cdot\f{\xi}{|\xi|}}g(\f{x}{|\xi|},\f{x}{|x|})
cos(\vec n, x_\ell)dS\\
&=\ds\f{\chi_\ell(\xi)}{|\xi|^n}\ds\f{1}{(i\f{\xi_\ell}{|\xi|})^m}\int_{1\le
|x|\le |\xi|}e^{-ix \cdot\f{\xi}{|\xi|}}\p_l^m
\biggl(g(\f{x}{|\xi|},\f{x}{|x|})\biggr)dx\\
&\quad +\ds\sum_{k=0}^{m-1}
\ds\f{\chi_\ell(\xi)}{|\xi|^n}\ds\f{1}{(i\f{\xi_\ell}{|\xi|})^k}
\f{1}{(-i\f{\xi_\ell}{|\xi|})}\int_{|x|=1}e^{-ix
\cdot\f{\xi}{|\xi|}}\p_\ell^k\biggl(g(\f{x}{|\xi|},\f{x}{|x|})\biggr)
cos(\vec n, x_\ell)dS\\
&=III+IV.\tag2.3
\endalign
$$

Due to
$$|\p_\ell^k\biggl(g(\ds\f{x}{|\xi|},\ds\f{x}{|x|})\biggr)|\le
\ds\sum_{j=0}^k\ds\sum_{|\al|\le k-j}C_{\al
j}|(\p_\ell^j\p_y^{\al}g)(\f{x}{|\xi|},\f{x}{|x|})| \ |x|^{-(k-j)} \
|\xi|^{-j},$$ then from (2.1), $IV$ is dominated by the multiplier
of $|\xi|^{-n}$, and moreover,
$$
\align
|III|&\le \ds\f{C}{|\xi|^n}\ds\sum_{\al+\beta=m}C_{\al\beta}\f{1}{|\xi|^{\al}}
\int_{1\le |x|\le |\xi|}|x|^{-\beta}dx\\
&\le\cases
C\ds\sum_{\al+\beta=m}C_{\al\beta}\biggl(\f{1}{|\xi|^m}+\f{1}{|\xi|^{n+\al}}\biggr)\quad\text
{if $\beta\not=n$};\\
\ds\f{C}{|\xi|^{n}}ln|\xi|\quad\text
{if $\beta=n$}.\\
\endcases\tag2.4
\endalign$$

Therefore, for $m\ge n$ and $|\xi|\ge 1$, we have
$|\hat{\vp}(\xi)|\le\ds\f{C(1+ln|\xi|)}{|\xi|^{n}}$ by (2.2)-(2.4),
which derives $(1+|\xi|)^{\f{n}{2}-\dl}\hat{\vp}(\xi)\in
L^2(\{|\xi|\ge 1 \})$ for any $\dl>0$, and further completes the
proof of (i).

(ii) Without loss of generality, we assume $supp \ \psi_i(x)\subset
[-1,1; -1,1]$ ($1\le i\le 4$).

Since
$$
\align |\hat{\vp}(\xi)|\equiv &
\biggl|\int_0^1\int_0^1\psi_1(x)e^{-ix\cdot\xi}dx
+\int_{-1}^0\int_0^1\psi_2(x)e^{-ix\cdot\xi}dx
\\
& \qquad \qquad +\int_{-1}^0\int_{-1}^0\psi_3(x)e^{-ix\cdot\xi}dx
+\int_0^1\int_{-1}^0\psi_4(x)e^{-ix\cdot\xi}dx\biggr|\\
\le & \cases
\ds\f{C}{|\xi_1| \ |\xi_2|}\quad\text{for $|\xi_1|\ge 1, |\xi_2|\ge 1$;}\\
\ds\f{C}{|\xi_1|}\quad\text{for $|\xi_1|\ge 1, |\xi_2|<1$;}\\
\ds\f{C}{|\xi_2|}\quad\text{for $|\xi_1|<1, |\xi_2|\ge 1$;}\\
C\quad\text{for $|\xi_1|<1, |\xi_2|<1$,}\\
\endcases\
\endalign$$
then from the fact  $1+|\xi| \le (1+|\xi_1|)(1+|\xi_2|) $ one has for
any $0<\dl<\f12$
$$
\align
&\int_{\Bbb R^2}(1+|\xi|)^{1-\dl}|\hat{\vp}(\xi)|^2d\xi\\
&\le
C\ds\prod_{i=1}^2\int_{1}^{\infty}\f{(1+|\xi_i|)^{1-\dl}}{|\xi_i|^2}d\xi_i
+C\ds\sum_{i=1}^2\int_{1}^{\infty}\f{(1+|\xi_i|)^{1-\dl}}{|\xi_i|^2}d\xi_i+C\\
&\le C.
\endalign
$$

Thus, the proof of (ii) is completed.

(iii) The proof procedure is similar to that in (ii), we omit it
here. \hfill \qed

\vskip .2cm

{\bf Remark 2.1.} {\it By the similar proof procedure as in Lemma
2.1.(i), we can also prove:  If $f(x)\in C^{\infty}(\Bbb
R^n\setminus\{0\})$ and has compact support, moreover,
$|\p^{\al}f(x)|\le C_{\al}|x|^{r-|\al|}$ for $x\not=0$ and
$r>-\f{n}{2}$, then $f(x)\in H^{\f{n}{2}+r-}(\Bbb R^n)$.}

\vskip .2cm

{\bf Remark 2.2.} {\it Under the assumption $(A_2)$, for any $\al\in
(\Bbb N\cup\{0\})^{n-1}$, we can also have that
$\p_{x'}^{\al}\vp(x)\in H^{\f12-\dl}(\Bbb R^n)$ for any $\dl >0$
small, here $x'=(x_2, ..., x_n)$. Thus,
$(1+|\xi|)^{\frac{1}{2}-\delta} (1+|\xi'|)^{|\alpha|} \hat{\vp}(\xi)
\in L^2(\Bbb R^n)$, where $\xi'=(\xi_2, \cdots, \xi_n)$.}

\vskip .2cm

 {\bf Lemma 2.2.} {\it If $u(t,x)\in C([0, T],
H^{\f12-}(\Bbb R^n))$ is a solution of the following linear equation
$$
\cases &\p_t^2u-t^m\Delta u=0,\qquad (t,x)\in [0,
+\infty)\times\Bbb R^n,\\
&u(0,x)=\psi(x),\quad \p_tu(0,x)=\vp(x),
\endcases\tag2.5
$$ where $\vp(x)$
satisfies the assumption $(A_2)$, $\p_{x'}^{\al}\psi(x)\in
H^{\f32-}(\Bbb R^n)$ for all $0\le |\al|\le [\f{n}{2}]+1$,  then
$u(t,x)\in L^{\infty}([0, T]\times\Bbb R^n)$.}

\vskip .1cm {\bf Proof.} Set $y(t,\xi)=\int_{\Bbb
R^n}u(t,x)e^{-ix\cdot\xi}dx$ with $\xi\in\Bbb R^n$ and
$y''(t,\xi)\equiv\p_t^2y(t,\xi)$, then it follows from the equation
of (2.5) that

$$y''(t,\xi)+t^{m}|\xi|^2y(t,\xi)=0.\tag2.6$$

Let $\tau=\ds\f{2t^{\f{m}{2}+1}|\xi|}{m+2}$ and $v(\tau)\equiv
y(t,|\xi|)$, then
$$\ds\f{d^2v}{d\tau^2}+\ds\f{m}{(m+2)\tau}\f{dv}{d\tau}+v=0.\tag2.7$$

As in [25], taking $z\equiv 2i\tau=\ds\f{4i}{m+2}t^{\f{m+2}{2}}|\xi|$ and
$w(z)=v(\ds\f{z}{2i})e^{\f{z}{2}}$ yields for $t>0$ and $|\xi|\neq
0$
$$zw''(z)+(\f{m}{m+2}-z)w'(z)-\f{m}{2(m+2)}w(z)=0.\tag2.8$$
(2.8) has two linearly independent solutions
$w_1(z)=\Phi(\f{m}{2(m+2)}, \f{m}{m+2}; z)$ and
$w_2(z)=z^{\f{2}{m+2}}\Phi(\f{m+4}{2(m+2)},\f{m+4}{m+2};z)$ by [13],
which are called the confluent hypergeometric functions.

By (2.6)-(2.8) and [23], we have for $t\ge 0$ and $\xi\in\Bbb R^n$
$$
\align y(t,\xi)&=
V_1(t,|\xi|)\psi^{\wedge}(\xi)+V_2(t,|\xi|)\vp^{\wedge}(\xi)\\
&\equiv y_1(t,\xi)+y_2(t,\xi)\tag2.9\endalign$$ with
$$
\cases &V_1(t,
|\xi|)=e^{-\f{z}{2}}\Phi(\f{m}{2(m+2)},\f{m}{m+2};z),\\
&V_2(t,|\xi|)=te^{-\f{z}{2}}\Phi(\f{m+4}{2(m+2)},\f{m+4}{m+2};z).
\endcases\tag2.10
$$

Since $\Phi(\f{m}{2(m+2)},\f{m}{m+2};z)$ and
$\Phi(\f{m+4}{2(m+2)},\f{m+4}{m+2};z)$ are analytic functions of
$z$, then $|\Phi(\f{m}{2(m+2)},\f{m}{m+2};z)|$ and
$|\Phi(\f{m+4}{2(m+2)},\f{m+4}{m+2};z)|\le C_M$ for $|z|\le M$. For
sufficiently large $|z|$, we have from formula (9) in pages 253 of
[13] that
$$|\Phi(\f{m}{2(m+2)},\f{m}{m+2};z)|\le C|z|^{-\f{m}{2(m+2)}}
\big(1+O\big(|z|^{-1}\big)\big),\quad
|\Phi(\f{m+4}{2(m+2)},\f{m+4}{m+2};z)|\le C|z|^{-\f{m+4}{2(m+2)}}
\big(1+O\big(|z|^{-1}\big)\big).\tag2.11$$

From Remark 2.2, we have that for $0\le |\al|\le [\f{n}{2}]+1$ and $0<\dl<\f12$
$$\hat{\vp}(\xi)=\ds\f{g_{\al}(\xi)}{(1+|\xi_1|)^{\f12-\dl}(1+|\xi'|)^{|\al|}},\tag2.12$$
where $g_{\al}(\xi)\in L^2(\Bbb R^n)$, $\xi'=(\xi_2, ...., \xi_n)$.

Therefore, for any $t\in (0, T]$, we have
$$
\align &\int_{\Bbb R^n} |y_2(t,\xi)|d\xi\le Ct\int_{\Bbb
R^n}\bigg|e^{-\f{2i}{m+2}t^{\f{m+2}{2}}|\xi|}
\Phi(\f{m+4}{2(m+2)},\f{m+4}{m+2};\f{4i}{m+2}t^{\f{m+2}{2}}|\xi|)
\vp^\wedge(\xi)\bigg|d\xi\\
&\quad \le Ct\big(\f{m+2}{2t^{\f{m+2}{2}}}\big)^{n} \int_{\Bbb
R^n}\bigg| \Phi(\f{m+4}{2(m+2)},\f{m+4}{m+2};2i|\eta|)
\vp^\wedge(\f{(m+2)\eta}{2t^{\f{m+2}{2}}})\bigg|d\eta\\
&\quad \le Ct\big(\f{m+2}{2t^{\f{m+2}{2}}}\big)^{n}
\int_{\Bbb R^n}\ds\f{1}{(1+|\eta|^2)^{\f{m+4}{4(m+2)}}}
|\vp^\wedge(\f{(m+2)\eta}{2t^{\f{m+2}{2}}})|d\eta\qquad\qquad\qquad \text{(by (2.11))}\\
&\quad \le
C_{\al}t\big(\f{m+2}{2t^{\f{m+2}{2}}}\big)^{\f{n}{2}}\biggl(\int_{\Bbb
R^n}\ds\f{1}{(1+|\eta|^2)^{\f{m+4}{2(m+2)}}}\ds\f{1}{(1+\f{|\eta_1|}{t^{\f{m+2}{2}}})^{1-2\dl}}
\ds\f{1}{(1+\f{|\eta'|}{t^{\f{m+2}{2}}})^{2|\al|}}
d\eta\biggr)^{\f12}\qquad \text{(by (2.12))}\\
&\quad \le C_{\al}t\biggl(\int_{\Bbb
R^n}\ds\f{1}{(1+t^{m+2}|\eta|^2)^{\f{m+4}{2(m+2)}}}\ds\f{1}{(1+|\eta_1|)^{1-2\dl}}
\ds\f{1}{(1+|\eta'|)^{2|\al|}}d\eta\biggr)^{\f12}\qquad \text{(choosing $|\al|=[\f{n}{2}]+1>\f{n}{2}$)}\\
&\quad \le C_Tt^{1-\f{3(m+2)\dl}{4}}\biggl(\int_{\Bbb
R}\f{1}{(1+|\eta_1|)^{1+\dl}}d\eta_1\biggr)^{\f12}\qquad
\text{(choosing $\dl<\ds\f{m+4}{3(m+2)}$)}\\
&\quad \le C_T \qquad\qquad \text{(choosing
$\dl<\ds\f{4}{3(m+2)}$)}.
\endalign
$$

Similarly,
$$
\align &\int_{\Bbb R^n} |y_1(t,\xi)|d\xi\le
C\big(\f{m+2}{2t^{\f{m+2}{2}}}\big)^{n} \int_{\Bbb R^n}\bigg|
\Phi(\f{m}{2(m+2)},\f{m}{m+2};2i|\eta|)
\psi^\wedge(\f{(m+2)\eta}{2t^{\f{m+2}{2}}})\bigg|d\eta\\
&\quad \le
C_{\al}\big(\f{m+2}{2t^{\f{m+2}{2}}}\big)^{\f{n}{2}}\biggl(\int_{\Bbb
R^n}\ds\f{1}{(1+|\eta|^2)^{\f{m}{2(m+2)}}}\ds\f{1}{(1+\f{|\eta_1|}{t^{\f{m+2}{2}}})^{3-2\dl}}
\ds\f{1}{(1+\f{|\eta'|}{t^{\f{m+2}{2}}})^{2|\al|}}
d\eta\biggr)^{\f12}\\
&\quad \le C_{\al}\biggl(\int_{\Bbb
R^n}\ds\f{1}{(1+|\eta_1|)^{3-2\dl}}
\ds\f{1}{(1+|\eta'|)^{2|\al|}}d\eta\biggr)^{\f12}\\
&\quad \le C_T\qquad \text{(choosing $|\al|=[\f{n}{2}]+1$ and
$0<\dl<1$)}.
\endalign
$$

Thus, $|u(t,x)|\le \int_{\Bbb R^n}|y(t,\xi)|d\xi\le \int_{\Bbb
R^n}|y_1(t,\xi)|d\xi+ \int_{\Bbb R^n}|y_2(t,\xi)|d\xi\le C_T$ for
$(t,x)\in (0, T]\times\Bbb R^n$, and then Lemma 2.2 is proved.
\hfill \qed

\vskip .2cm

 {\bf Lemma 2.3.} {\it If $f(t,x)\in C([0,T], H^{s}(\Bbb
R^n))$ and $\p_{x'}^{\al}f(t,x)\in L^{\infty}([0,T], H^s(\Bbb R^n))$
with $s>\f12$ and  $|\al|\le [\f{n}{2}]+1$, $v(t,x)$ is a solution
to the following problem
$$
\cases &\p_t^2u-t^m\Delta u=f(t,x),\\
&u(0,x)=\p_tu(0,x)=0,
\endcases\tag2.13
$$
then $u(t,x)\in L^{\infty}([0, T]\times\Bbb R^n)$.}

\vskip .1cm

{\bf Proof.} By the assumptions on $f(t,x)$, we have
$$f^{\wedge}(t,\xi)=\ds\f{g_{\al}(t,\xi)}{(1+|\xi_1|)^{s}(1+|\xi'|)^{|\al|}},$$
where $g_{\al}(t,\xi)\in L^{\infty}([0, T], L^2(\Bbb R^n))$ and
$|\al|=[\f{n}{2}]+1$.

From (2.13), we have
$$u(t,x)=\biggl(\int_0^t (V_2(t,|\xi|)V_1(\tau,|\xi|)
-V_1(t,|\xi|)V_2(\tau,|\xi|))f^\wedge(\tau,\xi)d\tau\biggr)^{\vee}(t,x),$$
where the expressions of $V_1(t,|\xi|)$ and $V_2(t,|\xi|)$ are given
in (2.10).

It is noted that
$$
\align |u^{\wedge}(t,\xi)|&\le\int_0^t|
V_2(t,|\xi|)V_1(\tau,|\xi|)f^\wedge(\tau,\xi)|d\tau+\int_0^t
|V_1(t,|\xi|)V_2(\tau,|\xi|)f^\wedge(\tau,\xi)|d\tau\\
\equiv& I+II.\tag2.14
\endalign
$$

Set $\eta=\ds\f{2}{m+2}t^{\f{m+2}{2}}\xi$,  we have
$$
\align &|I|\le Ct\int_0^t |\Phi(\f{m+4}{2(m+2)}, \f{m+4}{m+2};
2i|\eta|)||\Phi(\f{m}{2(m+2)}, \f{m}{m+2};
2i(\f{\tau}{t})^{\f{m+2}{2}}|\eta|)|
|f^\wedge(\tau,\f{(m+2)\eta}{2t^{\f{m+2}{2}}})|d\tau\\
&\le C_{\al}t\int_0^t
(1+|\eta|)^{-\f{m+4}{2(m+2)}}(1+(\f{\tau}{t})^{\f{m+2}{2}}|\eta|)^{-\f{m}{2(m+2)}}
\ds\f{|g_{\al}(\tau,\f{(m+2)\eta}{2t^{\f{m+2}{2}}})|}{(1+\f{|\eta_1|}{t^{\f{m+2}{2}}})^{s}
(1+\f{|\eta'|}{t^{\f{m+2}{2}}})^{|\al|}}d\tau
\endalign$$
and thus
$$
\align \int_{\Bbb R^n}|I|&d\xi\le
C_{\al}t^{1-\f{(m+2)n}{4}}\int_0^td\tau \biggl(\int_{\Bbb R^n}
\ds\f{1}{(1+|\eta|)^{\f{m+4}{m+2}}(1+(\f{\tau}{t})^{\f{m+2}{2}}|\eta|)^{\f{m}{m+2}}(1+\f{|\eta_1|}{t^{\f{m+2}{2}}})^{2s}
(1+\f{|\eta'|}{t^{\f{m+2}{2}}})^{2|\al|}}d\eta\biggr)^{\f12}\\
&\le C_{\al}t\int_0^td\tau \biggl(\int_{\Bbb R^n}
\ds\f{1}{(1+|\eta_1|)^{2s}
(1+|\eta'|)^{2|\al|}}d\eta\biggr)^{\f12}\\
&\le C_{\al}t^2.\tag2.15
\endalign
$$

On the other hand, due to
$$|II|\le C_{\al}\int_0^t
(1+|\eta|)^{-\f{m}{2(m+2)}}(1+(\f{\tau}{t})^{\f{m+2}{2}}|\eta|)^{-\f{m+4}{2(m+2)}}
\ds\f{\tau
|g_{\al}(\tau,\f{(m+2)\eta}{2t^{\f{m+2}{2}}})|}{(1+\f{|\eta_1|}{t^{\f{m+2}{2}}})^{s}
(1+\f{|\eta'|}{t^{\f{m+2}{2}}})^{|\al|}}d\tau,$$ then we can obtain
as in (2.15)
$$
\int_{\Bbb R^n}|II|d\xi\le C  t^2.\tag2.16
$$

Substituting (2.15) and (2.16) into (2.14) yields
$$\int_{\Bbb R^n}|u^{\wedge}(t,
\xi)|d\xi\le C t^2.$$

Consequently, $|u(t,x)|\le \int_{\Bbb R^n}|u^{\wedge}(t,\xi)|d\xi\le
C t^2$, and the proof on Lemma 2.3 is completed. \hfill \qed

\vskip .2cm

 Finally, we study the $L^{\infty}$ property of solution
to the 2-D linear problem (1.1) under the assumption $(A_3)$.

\vskip .2cm

 {\bf Lemma 2.4.} {\it  If $u(t,x)\in C([0, T],
H^{\f12-}(\Bbb R^2))$ is a solution of the following linear problem
$$ \cases &\p_t^2u-t^m\Delta u=0,\qquad
(t,x)\in [0, T]\times\Bbb R^2,\\
&u(0,x)=0,\quad \p_tu(0,x)=\vp(x),
\endcases\tag2.17
$$
where $\vp(x)$ satisfies the assumption $(A_3)$, then $u\in
L^{\infty}([0, T]\times\Bbb R^2)$.}

{\bf Remark 2.3.} {\it Due to $\vp(x)\in H^{\f12-}(\Bbb R^2)$ by
Lemma 2.1.(ii), then the optimal regularity of the solution $u(t,x)$
to (2.17) is $L^{\infty}([0, T], H^{\f12+\f{2}{m+2}-}(\Bbb R^2))$
(see Proposition 3.3 in $\S 3$ below). Thus, for $m\ge 2$, we can
not derive $u(t,x)\in L^{\infty}([0, T]\times\Bbb R^2)$ directly by
the Sobolev imbedding theorem. On the other hand, the proof
procedure on Lemma 2.4 will be rather useful in analyzing the
singularity structure of $u(t,x)$ in $\S 6$ below.}

\vskip .1cm

{\bf Proof.} In terms of Corollary 3.5 in [24], we have the
following expression for the solution of (2.17)
$$
u(t,x)=2tC_m(\phi(1))^{\phi(1)}F(\g,\g;1;1)\int_0^1(1-s^2)^{-\g}(\p_tv)(s\phi(t),x)ds,\tag2.18
$$
where $C_m=(\ds\f{2}{m+2})^{\f{m}{m+2}}2^{-\f{2}{m+2}}$,
$\g=\ds\f{m}{2(m+2)}$, $F(\g,\g;1;1)=F(\g,\g;1;z)|_{z=1}$ with
$F(\g,\g;1;z)$ a hypergeometric function, which satisfies
$z(1-z)\o''(z)+(1-(2\g+1)z)\o'(z)-\g^2\o(z)=0$, and $v(t,x)$ is a
solution to the following linear wave equation
$$\p_t^2v-\Delta v=0,\quad v(0,x)=0,\quad \p_tv(0,x)=\vp(x).\tag2.19$$

From (2.19),  we have
$$v(t,x)=\ds\f{1}{2\pi}\int_{B(x,t)}\ds\f{\vp(\xi)}{\sqrt{t^2-(x_1-\xi_1)^2-(x_2-\xi_2)^2}}d\xi.\tag2.20$$

Let $v_i(t,x)$ be the smooth solution to the linear wave equation
$\p_t^2v_i-\Delta v_i=0$ with the initial data $(v_i(0,x),$
$\p_tv_i(0,x))$ $=(0, \psi_i(x))$. Then it follows from (2.20) and a
direct computation that for $t>0$ and $x_1>0, x_2>0$ (in other
domains, the expressions are completely analogous)
$$
v(t,x)=\cases
v_1(t,x)\quad\text{for $\ds\f{x_1}{t}\ge 1$, $\ds\f{x_2}{t}\ge 1$};\\
v_1(t,x)+I_1(t,x)\quad\text{for $\ds\f{x_1}{t}\le 1$,
$\ds\f{x_2}{t}\ge 1$};\\
v_1(t,x)+I_2(t,x)\quad\text{for $\ds\f{x_2}{t}\le 1$, $\ds\f{x_1}{t}\ge 1$};\\
\ds v_1(t,x)+I_1(t,x)+I_2(t,x)\qquad\text{for $0<x_1<t$, $0<x_2<t$, $|x|>t$};\\
v_1(t,x)+I_1(t,x)+I_2(t,x)+I_3(t,x)
\quad\text{for $x_1>0$, $x_2>0$, $|x|<t$}\\
\endcases\tag2.21$$
with
$$
\cases
&I_1(t,x)=\ds\f{1}{2\pi}\int_{x_1}^trdr\int_{-arccos(\f{x_1}{r})}^{arccos(\f{x_1}{r})}\ds\f{(\psi_2-\psi_1)(x-r\o)}
{\sqrt{t^2-r^2}}d\th,\\
&I_2(t,x)=\ds\f{1}{2\pi}\int_{x_2}^trdr\int_{-arccos(\f{x_2}{r})}^{arccos(\f{x_2}{r})}\ds\f{(\psi_4-\psi_1)(x-r\o)}
{\sqrt{t^2-r^2}}d\th,\\
&I_3(t,x)=\ds\f{1}{2\pi}\int_{|x|}^trdr\int_{arcsin(\f{x_1}{r})}^{arccos(\f{x_2}{r})}\ds\f{(\psi_1+\psi_3-\psi_2-\psi_4)(x-r\o)}
{\sqrt{t^2-r^2}}d\th,
\endcases
$$
where $\o=(cos\th, sin\th)$, $r=\sqrt{|x_1-\xi_1|^2+|x_2-\xi_2|^2}$
and $(x_1-\xi_1, x_2-\xi_2)=(r cos\th, r sin\th)$.

Due to $\vp(x)\in H^{\f12-}(\Bbb R^2)$ by Lemma 2.1.(ii), then it
follows from the regularity theory of solution to linear wave
equation that
$$v(t,x)\in C([0, T], H^{\f32-}(\Bbb R^2))\cap C^1([0, T], H^{\f12-}(\Bbb R^2))\subset
W^{1,1}([0, T]\times\Bbb R^2).$$

Thus, we can take the first order derivative $\p_t v$ piecewisely
for $t>0$  and $x_1>0, x_2>0$ as follows
$$
\p_tv(t,x)=\cases
\p_tv_1(t,x)\quad\text{for $\ds\f{x_1}{t}\ge 1$, $\ds\f{x_2}{t}\ge 1$};\\
\p_tv_1(t,x)+\p_tI_1(t,x)\quad\text{for $\ds\f{x_1}{t}\le 1$,
$\ds\f{x_2}{t}\ge 1$};\\
\p_tv_1(t,x)+\p_tI_2(t,x)\quad\text{for $\ds\f{x_2}{t}\le 1$, $\ds\f{x_1}{t}\ge 1$};\\
\p_tv_1(t,x)+\p_tI_1(t,x)+\p_tI_2(t,x)\qquad\text{for $0<x_1<t$, $0<x_2<t$, $|x|>t$};\\
\p_tv_1(t,x)+\p_tI_1(t,x)+\p_tI_2(t,x)+\p_tI_3(t,x)
\quad\text{for $x_1>0$, $x_2>0$, $|x|<t$}.\\
\endcases\tag2.22$$

Here we only treat the term $\p_tI_3$ in (2.22) since the treatments
on $\p_tI_1$ and $\p_tI_2$ are analogous or even simpler in their
corresponding domains.

If we set $\psi=\psi_1+\psi_3-\psi_2-\psi_4$, then it follows from a
direct computation that for $x_1>0$, $x_2>0$ and $|x|<t$
$$
I_3(t,x)
=\int_{x_1-\sqrt{t^2-x_2^2}}^0d\xi_1\int_{x_2-\sqrt{t^2-(x_1-\xi_1)^2}}^0\f{\psi(\xi)}
{\sqrt{t^2-|x-\xi|^2}}d\xi_2.\tag2.23
$$

Taking the transformations $x=ty$ and $\xi=t\eta$ in (2.23) yields
$$I_3(t,ty)=t J(t, \f{x}{t})$$
where
$$J(t, z)=\int_{z_1-\sqrt{1-z_2^2}}^0d\eta_1\int_{z_2-\sqrt{1-(z_1-\eta_1)^2}}^0\f{\psi(t\eta)}
{\sqrt{1-|z-\eta|^2}}d\eta_2\qquad\text{for $0<|z|<1$ and $z_1,
z_2>0$},$$ which derives
$$\p_1 I_3(t,x)=(\p_1 J)(t, \f{x}{t}), \qquad \p_2 I_3(t,x)=(\p_2 J)(t, \f{x}{t})$$
and thus
$$\p_tI_3(t,x)=\f{I_3(t,x)}{t}+\f{1}{t}\int_{x_1-\sqrt{t^2-x_2^2}}^0d\xi_1
\int_{x_2-\sqrt{t^2-(x_1-\xi_1)^2}}^0\f{\xi\cdot\na_{\xi}\psi(\xi)}
{\sqrt{t^2-|x-\xi|^2}}d\xi_2-\ds\f{x\cdot \na_x
I(t,x)}{t}.\tag2.24$$

It is noted that for $x_1>0$, $x_2>0$ and $|x|<t$,
$$
\align |\p_1&I_3(t,x)|=|\ds\lim_{h\to
0}\ds\f{I_3(t,x_1+h,x_2)-I_3(t,x_1, x_2)}{h}| \\= &
\biggl|\int_{x_1-\sqrt{t^2-x_2^2}}^0d\xi_1\int_{x_2-\sqrt{t^2-(x_1-\xi_1)^2}}^0\f{\p_1\psi(\xi)}
{\sqrt{t^2-|x-\xi|^2}}d\xi_2 -\int_{x_2-\sqrt{t^2-x_1^2}}^0\f{\psi(0,\xi_2)}{\sqrt{t^2-x_1^2-(x_2-\xi_2)^2}}d\xi_2\biggr|\\
= &
\biggl|\int_{x_1-\sqrt{t^2-x_2^2}}^0d\xi_1\int_{x_2-\sqrt{t^2-(x_1-\xi_1)^2}}^0\f{\p_1\psi(\xi)}
{\sqrt{t^2-|x-\xi|^2}}d\xi_2 +\int_{1}^{\f{x_2}{\sqrt{t^2 -x_1^2}}}\psi(0,x_2 -s\sqrt{t^2-x_1^2})d(arcsin s)\biggr|\\
=&\biggl|\int_{x_1-\sqrt{t^2-x_2^2}}^0d\xi_1\int_{x_2-\sqrt{t^2-(x_1-\xi_1)^2}}^0\f{\p_1\psi(\xi)}
{\sqrt{t^2-|x-\xi|^2}}d\xi_2+\psi(0, 0)arcsin(\ds\f{x_2}{\sqrt{t^2-x_1^2}})\\
&\qquad-\f{\pi}{2}\psi(0,x_2-\sqrt{t^2-x_1^2})+\sqrt{t^2-x_1^2}\int_1^{\f{x_2}{\sqrt{t^2-x_1^2}}}
\p_2\psi(0,x_2-s\sqrt{t^2-x_1^2})arcsin s\ ds\biggr|\\
\le &
C_T\biggl(1+\int_{x_1-\sqrt{t^2-x_2^2}}^0d\xi_1\int_{x_2-\sqrt{t^2-(x_1-\xi_1)^2}}^0\f{1}
{\sqrt{t^2-|x-\xi|^2}}d\xi_2\biggr)\\
\le & C_T(1+t)\tag2.25
\endalign$$

and
$$
\align |\p_2I_3(t,x)|=&
\biggl|\int_{x_1-\sqrt{t^2-x_2^2}}^0d\xi_1\int_{x_2-\sqrt{t^2-(x_1-\xi_1)^2}}^0\f{\p_2\psi(\xi)}
{\sqrt{t^2-|x-\xi|^2}}d\xi_2+\psi(0,0)arcsin(\ds\f{x_1}{\sqrt{t^2-x_2^2}})\\
&\qquad -\f{\pi}{2}\psi(x_1-\sqrt{t^2-x_2^2},0)
+\sqrt{t^2-x_2^2}\int_1^{\f{x_1}{\sqrt{t^2-x_2^2}}}\p_1\psi(x_1-s\sqrt{t^2-x_2^2},0)arcsin s \ ds\biggr|\\
\le & C_T(1+t)\tag2.26
\endalign
$$

On the other hand, analogous computation yields for $x_1>0$, $x_2>0$
and $|x|<t\le T$
$$|\f{I_3(t,x)}{t}|\le C_T\quad\text{and \quad $\biggl|\f{1}{t}\int_{x_1-\sqrt{t^2-x_2^2}}^0d\xi_1
\int_{x_2-\sqrt{t^2-(x_1-\xi_1)^2}}^0\f{\xi\cdot\na_{\xi}\psi(\xi)}
{\sqrt{t^2-|x-\xi|^2}}d\xi_2\biggr|\le C_T$}.\tag2.27$$

Therefore, $\p_tI_3(t,x)\in L^{\infty}$ in the domain $\{(t,x):
x_1>0, x_2>0, |x|<t\le T\}$ by (2.24). Similarly, we can obtain
$\p_tI_1(t,x)\in L^{\infty}$ and $\p_tI_2(t,x)\in L^{\infty}$ in the
related domains, and thus $\p_tv(t,x)\in L^{\infty}([0, T]\times\Bbb
R^2)$. These, together with (2.18), yield
$$u(t,x)\in L^{\infty}([0,T]\times\Bbb R^2)\quad \text{and\quad
$\|u(t,\cdot)\|_{L^{\infty}(\Bbb R^2)}\le
Ct\dsize\sum_{i=1}^4\|\psi_i(x)\|_{C^1}$}.\tag2.28$$

Consequently, we complete the proof of Lemma 2.4. \hfill \qed

\vskip .2cm

{\bf Remark 2.4.} {\it It is not difficult that  by the expression
(2.21) of $v(t,x)$, one can get $v(t,x)\in
C^{\infty}((0,T]\times\Bbb
R^2\setminus\Sigma_0\cup\Sigma_1^{\pm}\cup\Sigma_2^{\pm})$, where
$\Sigma_0=\{(t,x): t>0, |x|=t\}$, $\Sigma_1^{\pm}=\{(t,x): t>0,
x_1=\pm t\}$ and $\Sigma_2^{\pm}=\{(t,x): t>0, x_2=\pm t\}$. On the other hand, $v(t,x)\not\in
C^2((0,T]\times\Bbb R^2)$ since $v(t,x)$ has a strong singularity
when the variables $(t,x)$ go across
$\Sigma_0\cup\Sigma_1^{\pm}\cup\Sigma_2^{\pm}$. Indeed, for example,
it follows from (2.25) and a direct computation that for $x_1>0,
x_2>0$ and $|x|<t$
$$
\align
\p_{12}^2I_3&(t,x)=\int_{x_1-\sqrt{t^2-x_2^2}}^0d\xi_1\int_{x_2-\sqrt{t^2-(x_1-\xi_1)^2}}^0\f{\p_{12}^2\psi(\xi)}
{\sqrt{t^2-|x-\xi|^2}}d\xi_2+\ds\f{3\psi(0,0)}{\sqrt{t^2-|x|^2}}\\
&\qquad
 -\f{\pi}{2}\p_1\Big(\psi(x_1-\sqrt{t^2-x_2^2},0) \Big)+\sqrt{t^2-x_2^2}\int_1^{\f{x_1}{\sqrt{t^2-x_2^2}}}\p_1^2\psi(x_1-s\sqrt{t^2-x_2^2},0)arcsin
sds\\
&\qquad -\f{\pi}{2}\p_2\psi(0,x_2-\sqrt{t^2-x_1^2})
+\sqrt{t^2-x_1^2}\int_1^{\f{x_2}{\sqrt{t^2-x_1^2}}}\p_2^2\psi(0,x_2-s\sqrt{t^2-x_1^2})arcsin
sds \\ &=\ds\f{3\psi(0,0)}{\sqrt{t^2-|x|^2}}+\text{\bf bounded
terms}.\tag2.29
\endalign
$$
Thus, (2.29) implies $\p_{12}^2I_3(t,x)\to\infty$ as
$(t,x)\to\Sigma_0$ since $\psi(0)\not =0$ can be assumed without
loss of generality (this is due to the assumption of
$\psi_i(0)\not=\psi_j(0)$ for some $i\not=j$ and $1\le i<j\le 4$ in
$(A_3)$ and the different expressions of $\psi(x)$ in the related
domains $\{(t,x): t>0, \pm x_1>0, \pm x_2>0\}$ respectively). In
addition, by an analogous computation, we can derive that
$\p_{12}^2I_1(t,x)$ and $\p_{12}^2I_2(t,x)$ are bounded for $x_1>0,
x_2>0$ and $|x|<t$. Hence $\p_{12}^2v(t,x)\to\infty$ as
$(t,x)\to\Sigma_0$ and further $v(t,x)\not\in C^2((0,T]\times\Bbb
R^2)$ is proved. However, by the expression (2.18) and due to the
lack of strong Huyghens' principle for the Tricomi-type equations,
we can show that the solution $u(t,x)\not\in C^2((0, T]\times\Bbb
R^2\setminus\G_0\cup\G_1^{\pm}\cup\G_2^{\pm}) $ of (2.17) holds true
in $\S 6$ below, which implies an essential difference between the
degenerate equation and the strict hyperbolic equation.}

\vskip 0.3 true cm \centerline{\bf $\S 3$. Some regularity estimates
on the solutions to linear generalized Tricomi equations} \vskip 0.3
true cm

At first, we list some results on the  confluent hypergeometric
functions for our computations later on.

The confluent hypergeometric equation is
$$zw''(z)+(c-z)w'(z)-aw(z)=0,\tag3.1$$
where $z\in{\Bbb C}$, $a$ and $c$ are constants. The solution of
(3.1) is called the confluent hypergeometric function.

When $c$ {\bf is not an integer}, (3.1) has two {\bf linearly
independent} solutions:

$$w_1(z)=\Phi(a,c;z),\q
w_2(z)=z^{1-c}\Phi(a-c+1,2-c;z).$$

Below are some crucial properties of the confluent hypergeometric
functions.

\vskip .2cm

{\bf Lemma 3.1.}  {\it

(i) (pages 278 of [13])). For $-\pi<arg z<\pi$ and large $|z|$, then
$$\align \Phi(a,c;z)=&\f{\Gamma(c)}{\Gamma(c-a)}(e^{i\pi\epsilon}z^{-1})^a
\sum_{n=0}^M\f{(a)_n(a-c+1)_n}{n!}(-z)^{-n}
+O\big(|z|^{-a-M-1}\big)\\
\quad &+
\f{\Gamma(c)}{\Gamma(a)}e^{z}z^{a-c}\sum_{n=0}^N\f{(c-a)_n(1-a)_n}{n!}z^{-n}
+O\big(|e^zz^{a-c-N-1}|\big),\tag3.2\endalign$$ where $\epsilon=1$
if $Im z>0$, $\epsilon=-1$ if $Im z<0$, $(a)_0\equiv 1, (a)_n\equiv
a(a+1)\cdot\cdot\cdot (a+n-1)$, and $M,N=0,1,2,3...$.

(ii) (page 253 of [13]). $\Phi(a,c;z)=e^z\Phi(c-a,c;-z)$.

(iii) (page 254 of [13]).
$$\f{d^n}{dz^n}\Phi(a,c;z)=\f{(a)_n}{(c)_n}\Phi(a+n,c+n;z)\tag3.3$$
and
$$\f{d}{dz}\Phi(a,c;z)=\f{1-c}{z}\bigg(\Phi(a,c;z)-\Phi(a,c-1;z)\bigg).\tag3.4$$}

For such a problem
$$
\cases &\p_t^2u-t^m\Delta u=0,\qquad (t,x)\in [0,
+\infty)\times\Bbb R^n,\\
&u(0,x)=\phi_1(x),\quad \p_tu(0,x)=\phi_2(x),
\endcases\tag3.5
$$
by the results in [23], one has for $t\ge 0$
$$u^{\wedge}(t,\xi)=V_1(t,|\xi|)\phi_1^{\wedge}(\xi)+V_2(t,|\xi|)\phi_2^{\wedge}(\xi),\tag3.6$$
where the expressions of $V_1(t,|\xi|)$ and $V_2(t,|\xi|)$ have been
given in (2.10).

In order to analyze the regularities of $u^{\wedge}(t,\xi)$ in (3.6)
under some restrictions on $\phi_i(x) (i=1,2)$,  we require to
establish the following estimates:

\vskip .2cm

{\bf Lemma 3.2.} {\it For $0\le s_1\leq \f{m}{2(m+2)}$, $0\le
s_2\leq \f{m+4}{2(m+2)}$ and some fixed positive constant $T$, if
$g(x)\in H^s(\Bbb R^n)$ with $s\in\Bbb R$, then we have for $0<t\le
T$

(i) $\cases \|\big(V_1(t,|\xi|)
g^\wedge(\xi)\big)^\vee\|_{H^{s+s_1}}\leq C
t^{-\f{s_1(m+2)}{2}}\|g\|_{H^{s}},\\
\|\big(V_2(t,|\xi|) g^\wedge(\xi)\big)^\vee\|_{H^{s+s_2}}\leq C
t^{1-\f{s_2(m+2)}{2}}\|g\|_{H^{s}}.\\
\endcases\ \q\q\q\q\q\q\q\q\q\q\q\q\q\q\q\q\q\q\q\q(3.7)
$

(ii) $\cases \|\big(\p_t V_1(t,|\xi|)
g^\wedge(\xi)\big)^\vee\|_{H^{s-\f{m+4}{2(m+2)}}}\leq
C\|g\|_{H^{s}},\\
\|\big(\p_t V_2(t,|\xi|)
g^\wedge(\xi)\big)^\vee\|_{H^{s-\f{m}{2(m+2)}}}\leq C\|g\|_{H^{s}}
\endcases\ \q\q\q\q\q\q\q\q\q\q\q\q\q\q\q\q\q\q\q\q\q(3.8)
$}

\vskip .1cm

{\bf Proof.}  (i) First, we fix $t=(\f{m+2}{2})^{\f{2}{m+2}}$ to
show (3.7) (in this case, the corresponding variable $z$ in (2.10)
becomes $z=2i|\xi|$). Subsequently, for the variable $t$, as in [25]
and so on, we can use the scaling technique to derive (3.7).

Since $\Phi(a,c;z)$ is an analytic function of $z$, then
$\Phi(\f{m}{2(m+2)},\f{m}{m+2};2i|\xi|)$ and
$\Phi(\f{m+4}{2(m+2)},\f{m+4}{m+2};2i|\xi|)$ are bounded for
$|\xi|\leq C$. On the other hand, it follows from  (3.2) that for
large $|\xi|$

$$|\Phi(\f{m}{2(m+2)},\f{m}{m+2};2i|\xi|)|\leq
C(1+|\xi|^2)^{-\f{m}{4(m+2)}}$$ and

$$ |\Phi(\f{m+4}{2(m+2)},\f{m+4}{m+2};2i|\xi|)|\leq
C(1+|\xi|^2)^{-\f{m+4}{4(m+2)}}.$$

Thus, for any $s_1\in [0, \f{m}{2(m+2)}]$ and $s_2\in
[0,\f{m+4}{2(m+2)}]$, by a direct computation, we arrive at
$$\align &\|(V_1((\f{m+2}{2})^{\f{2}{m+2}},|\xi|)
g^\wedge(\xi))^\vee\|_{H^{s+s_1}}\\
=&\|(1+|\xi|^2)^{\f{s+s_1}{2}}e^{-i|\xi|}\Phi(\f{m}{2(m+2)},\f{m}{m+2};2i|\xi|)
g^\wedge(\xi)\|_{L^2}\\
\leq
&\|(1+|\xi|^2)^{\f{s_1}{2}}\Phi(\f{m}{2(m+2)},\f{m}{m+2};2i|\xi|)\|_{L^{\infty}}
\|(1+|\xi|^2)^\f{s}{2}g^\wedge(\xi)\|_{L^2}\\
\leq &C\|g\|_{H^s}\tag3.9\endalign$$ and

$$\|(V_2((\f{m+2}{2})^{\f{2}{m+2}},|\xi|)
g^\wedge(\xi))^\vee\|_{H^{s+s_2}}\leq C\|g\|_{H^s}.\tag3.10$$

Next we treat $\|(V_1(t,|\xi|) g^\wedge(\xi))^\vee\|_{H^{s+s_1}}$
and $\|(V_2(t,|\xi|) g^\wedge(\xi))^\vee\|_{H^{s+s_2}}$. To this
end, we introduce the following transformation
$$
\eta=\f{2}{m+2}t^{\f{m+2}{2}}\xi,$$ and then we have
$$\align &\|(V_1(t,|\xi|)
g^\wedge(\xi))^\vee\|_{H^{s+s_1}}\\
=&\biggl(\int_{\Bbb
R^n}\bigg|(1+|\xi|^2)^{\f{s_1}{2}}e^{-\f{2i}{m+2}t^{\f{m+2}{2}}|\xi|}
\Phi(\f{m}{2(m+2)},\f{m}{m+2};\f{4i}{m+2}t^{\f{m+2}{2}}|\xi|)
(1+|\xi|^2)^{\f{s}{2}}g^\wedge(\xi)\bigg|^2d\xi\biggr)^{\f{1}{2}}\\
=&\big(\f{m+2}{2t^{\f{m+2}{2}}}\big)^{\f{n}{2}}
\bigg(\int_{R^n}\bigg|(1+|\f{(m+2)\eta}{t^{\f{m+2}{2}}}|^2)^{\f{s_1}{2}}
\Phi(\f{m}{2(m+2)},\f{m}{m+2};2i|\eta|)
G^\wedge(\eta)\bigg|^2d\eta\bigg)^{\f{1}{2}}\tag3.11\endalign$$ and
$$\align &\quad\|(V_2(t,|\xi|)
g^\wedge(\xi))^\vee\|_{H^{s+s_2}}\\
=&t\big(\f{m+2}{2t^{\f{m+2}{2}}}\big)^{\f{n}{2}}
\bigg(\int_{R^n}\bigg|(1+|\f{(m+2)\eta}{2t^{\f{m+2}{2}}}|^2)^{\f{s_2}{2}}
\Phi(\f{m+4}{2(m+2)},\f{m+4}{m+2};2i|\eta|)
G^\wedge(\eta)\bigg|^2d\eta\bigg)^{\f{1}{2}},\tag3.12\endalign$$
here and below the notation $G^\wedge(\eta)$ is defined as

$$G^\wedge(\eta)=\biggl(1+\biggl|\f{(m+2)\eta}{2t^{\f{m+2}{2}}}\biggr|^2\biggr)^{\f{s}{2}}
g^\wedge(\f{(m+2)\eta}{2t^{\f{m+2}{2}}}).$$

It is noted that
$$
\|G^\wedge(\eta)\|_{L^2}
=(\int_{R^n}|(1+|\xi|^2)^{\f{s}{2}}g^\wedge(\xi)|^2
(\f{2t^{\f{m+2}{2}}}{m+2})^nd\xi)^{\f{1}{2}} \leq
Ct^{\f{n(m+2)}{4}}\|g\|_{H^s}.\tag3.13$$

Additionally, for $0<t\leq T$ and $\alpha\geq 0$, we have

$$\bigg(1+|\f{(m+2)\eta}{2t^{\f{m+2}{2}}}|^2\bigg)^{\alpha}
<Ct^{-\alpha(m+2)}(1+|\eta|^2)^{\alpha}.\tag3.14$$

Thus, we obtain from (3.11)-(3.14) that for $0<t\leq T$
$$\align &\|(V_1(t,|\xi|)
g^\wedge(\xi))^\vee\|_{H^{s+s_1}}\\
\leq&Ct^{-\f{s_1(m+2)}{2}-\f{n(m+2)}{4}}
(\int_{R^n}|(1+|\eta|^2)^{\f{s_1}{2}}\Phi(\f{m}{2(m+2)},\f{m}{m+2};2i|\eta|)
G^\wedge(\eta)|^2d\eta)^{\f{1}{2}}\\
=&Ct^{-\f{s_1(m+2)}{2}-\f{n(m+2)}{4}}\|
G^\wedge(\eta)\|_{L^2}\\
\leq&Ct^{-\f{s_1(m+2)}{2}}\|g\|_{H^s}\tag3.15\endalign$$ and
$$\|(V_2(t,|\xi|)
g^\wedge(\xi))^\vee\|_{H^{s+s_2}} \leq
Ct^{1-\f{s_2(m+2)}{2}}\|g\|_{H^s}.\tag3.16$$

Consequently, we complete the proof of Lemma 3.2.(i).

(ii). It follows from a direct computation and (3.3)-(3.4) that
$$\align &\p_t V_1(t,|\xi|)\\
=&2i(\f{m+2}{4i})^{\f{m}{m+2}}|\xi|^{\f{2}{m+2}}z^{\f{m}{m+2}}\bigg(-\f{1}{2}e^{-\f{z}{2}}
\Phi(\f{m}{2(m+2)},\f{m}{m+2},z)+\f{1}{2}e^{-\f{z}{2}}
\Phi(\f{3m+4}{2(m+2)},\f{2(m+1)}{m+2},z)\bigg)\\
=&i(\f{m+2}{4i})^{\f{m}{m+2}}|\xi|^{\f{2}{m+2}}z^{\f{m}{m+2}}
e^{-\f{z}{2}}\bigg(\Phi(\f{3m+4}{2(m+2)},\f{2(m+1)}{m+2},z)
-\Phi(\f{m}{2(m+2)},\f{m}{m+2},z)\bigg)\tag3.17\endalign$$ and
$$\p_t V_2(t,|\xi|)
=e^{-\f{z}{2}}\bigg(\Phi(\f{m+4}{2(m+2)},\f{2}{m+2},z)
-\f{(m+2)z}{4}\Phi(\f{m+4}{2(m+2)},\f{m+4}{m+2},z)\bigg).
\tag3.18$$

Thus, in terms of (3.2), we have for large $|z|$
$$\align &|\p_t V_1(t,|\xi|)|\\
\leq&
\big|i(\f{m+2}{4i})^{\f{m}{m+2}}|\xi|^{\f{2}{m+2}}z^{\f{m}{m+2}}
e^{-\f{z}{2}}\big|\bigg(|\Phi(\f{3m+4}{2(m+2)},\f{2(m+1)}{m+2},z)|
+|\Phi(\f{m}{2(m+2)},\f{m}{m+2},z)|\bigg)\\
\leq
&C|\xi|^{\f{2}{m+2}}|z|^{\f{m}{m+2}}\bigg(|z|^{max(-\f{m}{2(m+2)},-\f{3m+4}{2(m+2)})}
+|z|^{max(-\f{m}{2(m+2)},-\f{m}{2(m+2)})}\bigg)\\
\le &Ct^{\f{m}{4}}(1+|\xi|^2)^{\f{m+4}{4(m+2)}}\tag3.19\endalign$$
and

$$\align &|\p_t V_2(t,|\xi|)|
\leq |e^{-\f{z}{2}}|\bigg(|\Phi(\f{m+4}{2(m+2)},\f{2}{m+2},z)|
+|\f{(m+2)z}{4}\Phi(\f{m+4}{2(m+2)},\f{m+4}{m+2},z)|\bigg)\\
\leq& C\bigg(|z|^{max(\f{m}{2(m+2)},-\f{m+4}{2(m+1)})}
+|z|^{1+max(-\f{m+4}{2(m+2)},-\f{m+4}{2(m+2)})}\bigg)\\
\le& Ct^{\f{m}{4}}(1+|\xi|^2)^{\f{m}{4(m+2)}}.\tag3.20\endalign$$

Next it suffices to estimate $\|(\p_t V_2(t,|\xi|)
g^\wedge(\xi))^\vee\|_{H^{s-\f{m}{2(m+2)}}}$ since the treatment on
$(\p_t V_1(t,|\xi|) g^\wedge(\xi))^\vee$ is completely analogous.

As in (i), we fix $t=(\f{m+2}{2})^{\f{2}{m+2}}$. In this case, by
the analytic property of $\Phi(a,c ; z)$ and $(3.20)$, we arrive at

$$\align &\|(1+|\xi|^2)^{\f{s}{2}-\f{m}{4(m+2)}}\p_tV_2((\f{m+2}{2})^{\f{2}{m+2}},|\xi|)
g^\wedge(\xi)\|_{L^2}\\
\leq&\|(1+|\xi|^2)^{-\f{m}{4(m+2)}}\p_tV_2((\f{m+2}{2})^{\f{2}{m+2}},
|\xi|)\|_{L^{\infty}}\|(1+|\xi|^2)^\f{s}{2}g^\wedge\|_{L^2}\\
\leq &C\|g\|_{H^s}.\tag3.21\endalign$$

For any $t>0$, we have
$$\align \quad &\|(\p_t V_2(t,|\xi|)
g^\wedge(\xi))^\vee\|_{H^{s-\f{m}{2(m+2)}}}\\
=&\bigg(\int_{\Bbb
R^n}\bigg|(1+|\xi|^2)^{-\f{m}{4(m+2)}}e^{-\f{2i}{m+2}t^{\f{m+2}{2}}|\xi|}
\bigg(\Phi(\f{m+4}{2(m+2)},\f{2}{m+2},\f{4i}{m+2}t^{\f{m+2}{2}}|\xi|)\\
\quad
&-2it^{k+1}|\xi|\Phi(\f{m+4}{2(m+2)},\f{m+4}{m+2},\f{4i}{m+2}t^{\f{m+2}{2}}|\xi|)\bigg)
(1+|\xi|^2)^{\f{s}{2}}g^\wedge(\xi)\bigg|^2d\xi\bigg)^{\f{1}{2}}\\
=&(\f{m+2}{2t^{\f{m+2}{2}}})^{\f{n}{2}}
\bigg(\int_{R^n}\bigg|(1+|\f{(m+2)\eta}{2t^{\f{m+2}{2}}}|^2)^{-\f{m}{4(m+2)}}
\bigg(\Phi(\f{m+4}{2(m+2)},\f{2}{m+2},2i|\eta|)\\
\quad
&-i(m+2)|\eta|\Phi(\f{m+4}{2(m+2)},\f{m+4}{m+2},2i|\eta|)\bigg)
G^\wedge(\eta)\bigg|^2 d\eta\bigg)^{\f{1}{2}}.\tag3.22\endalign$$

Substituting the estimates $(3.13)$ and $(3.21)$ into the expression
$(3.22)$ yields for $0<t\le T$

$$\align &\|(\p_t V_2(t,|\xi|)
g^\wedge(\xi))^\vee\|_{H^{s-\f{m}{2(m+2)}}}\\
=&Ct^{-\f{n(m+2)}{4}}
\|(1+|\f{(m+2)\eta}{2t^{\f{m+2}{2}}}|^2)^{-\f{m}{4(m+2)}}
\p_tV_2((\f{m+2}{2})^{\f{2}{m+2}},|\eta|)G^\wedge(\eta)\|_{L^2}\\
\leq&Ct^{-\f{n(m+2)}{4}}
\|(1+|\eta|^2)^{-\f{m}{4(m+2)}}\p_tV_2((\f{m+2}{2})^{\f{2}{m+2}},|\eta|)\|_{L^{\infty}}
\|G^\wedge(\eta)\|_{L^2}\\
\leq&Ct^{-\f{n(m+2)}{4}}\|
G^\wedge(\eta)\|_{L^2}\\
\leq&C\|g\|_{H^s}.\endalign$$

By a similar way, the following estimate of $\|(\p_t V_1(t,|\xi|)
g^\wedge(\xi))^\vee\|_{H^{s-\f{m+4}{2(m+2)}}}$ can be obtained for
$0<t\leq T$
$$ \|(\p_t V_1(t,|\xi|)
g^\wedge(\xi))^\vee\|_{H^{s-\f{m+4}{2(m+2)}}} \leq C\|g\|_{H^s}.$$

Thus, we complete the proof of Lemma 3.2.(ii). \hfill \qed

\vskip .2cm

Based on Lemma 3.2, we can derive the following estimates for the
solution of (3.5).

\vskip .2cm

{\bf Proposition 3.3.} {\it If $\phi_1\in H^{s}(\Bbb R^n)$ and
$\phi_2\in H^{s-\f{2}{m+2}}(\Bbb R^n)$ with $s\in\Bbb R$, then (3.5)
has a solution $u(t,x)\in C([0, T], H^s(\Bbb R^n ))\cap C((0,T],
H^{s+\f{m}{2(m+2)}}(\Bbb R^n))\cap C^1([0, T],
H^{s-\f{m+4}{2(m+2)}}(\Bbb R^n))$ which admits the following
estimates for $0<t\leq T$
$$
\|u(t,\cdot)\|_{H^{s}(\Bbb
R^n)}+t^{\f{m}{4}}\|u(t,\cdot)\|_{H^{s+\f{m}{2(m+2)}}(\Bbb R^n)}
+\|\p_tu(t,\cdot)\|_{H^{s-\f{m+4}{2(m+2)}}(\Bbb R^n)}\leq
C(\|\phi_1\|_{H^s(\Bbb R^n)}+\|\phi_2\|_{H^{s-\f{2}{m+2}}(\Bbb
R^n)}).$$}

{\bf Proof.}  In Lemma 3.2, if we take $s_1=0$ or
$s_1=\f{m}{2(m+2)}$ and $s_2=0$ or $s_2=\f{m+4}{2(m+2)}$
respectively, then Proposition 3.3 can be shown. \hfill \qed

\vskip .2cm

Next, we consider the following inhomogeneous problem
$$\cases \p_t^2u-t^{m}\triangle u=f(t,x)\\
u(0,x)=0,\quad u_t(0,x)=0.\endcases\tag3.23$$

As in (2.14) of Lemma 2.3, one has
$$u^{\wedge}(t,\xi)=\int_0^t (V_2(t,|\xi|)V_1(\tau,|\xi|)
-V_1(t,|\xi|)V_2(\tau,|\xi|))f^\wedge(\tau,\xi)d\tau.\tag3.24$$

Based on Lemma 3.1-Lemma 3.2,  we can establish

\vskip .2cm

 {\bf Lemma 3.4.} {\it If $f(t,\cdot)\in C([0,T],
H^s(\Bbb R^n))$ with $s\in\Bbb R$ and $T$ a fixed positive constant,
then for $t\in (0, T]$
$$ \cases &\|u(t,\cdot)\|_{H^{s+p_1}}\leq Ct^{2-\f{
p_1}{2}(m+2)}\|f(t,x)\|_{L^{\infty}([0, T], H^s)},\quad\quad
\quad\quad \quad\quad \quad\quad \quad\quad \quad\quad \quad\quad
\quad\quad\quad\quad\quad\quad (3.25)\\
&\|\p_tu(t,\cdot)\|_{H^{s-\f{m}{2(m+2)}+p_2}}\leq
C_{p_2}t^{1-\f{m+2}{2}p_2}\|f(t,x)\|_{L^{\infty}([0, T],
H^s)},\quad\quad \quad\quad \quad\quad \quad\quad \quad\quad
\quad\quad \quad\quad\quad\quad(3.26)\\
\endcases
$$
where $0\le p_1<p_1(m)=\cases \ds\f{m+8}{2(m+2)}\quad\text{for $m\ge 4$}\\
1\quad\text{for $m\le 4$}\endcases$ and
$p_2<p_2(m)=\min\{\f{2}{m+2}, \f{m}{2(m+2)}\}$.}

{\bf Proof.} It follows from the Minkowski inequality and (3.24)
that
$$\align \quad \|u(t,\cdot)\|_{H^{s+p_1}}
\leq &\int_0^t\biggl(\int_{\Bbb
R^n}|(1+|\xi|^2)^{\f{s}{2}+\f{p_1}{2}} (V_2(t,|\xi|)V_1(\tau,|\xi|)
-V_1(t,|\xi|)V_2(\tau,|\xi|))f^\wedge(\tau,\xi)|^2d\xi\biggr)^{\f{1}{2}}d\tau\\
\leq&\int_0^t\|(1+|\xi|^2)^{\f{s}{2}+\f{p_1}{2}}
V_2(t,|\xi|)V_1(\tau,|\xi|)f^\wedge(\tau,\xi)\|_{L^2}d\tau\\
\quad&+\int_0^t\|(1+|\xi|^2)^{\f{s}{2}+\f{p_1}{2}}
V_1(t,|\xi|)V_2(\tau,|\xi|)f^\wedge(\tau,\xi)\|_{L^2}d\tau\\
\equiv&I_1+I_2.\tag3.27\endalign$$

Let $p_1=s_1+s_2$ with $0\le s_1<\min\{\f{m}{2(m+2)}, \f{2}{m+2}\}$
and $0\le s_2\le\f{m+4}{2(m+2)}$, then we have by Lemma 3.2
$$
\align I_1&\le
Ct^{1-\f{s_2(m+2)}{2}}\int_0^t\|(1+|\xi|^2)^{\f{s_1}{2}}V_1(\tau,|\xi|)(1+|\xi|^2)^{\f{s}{2}}f^{\wedge}(\tau,
\xi)\|_{L^2}d\tau\\
&\le
Ct^{1-\f{s_2(m+2)}{2}}\int_0^t\tau^{-\f{s_1(m+2)}{2}}\|f(\tau,\cdot)\|_{H^s}d\tau\\
&\le Ct^{2-\f{p_1(m+2)}{2}}\|f\|_{L^{\infty}([0, T],H^s)}.\tag3.28
\endalign
$$

On the other hand, if we set $p_1=\t s_1+\t s_2$ with $0\le \t
s_1\le\f{m}{2(m+2)}$ and $0\le \t s_2<\min\{\f{m+4}{2(m+2)},
\f{4}{m+2}\}$, then we have by Lemma 3.2
$$
\align I_2&\le Ct^{-\f{\t
s_1(m+2)}{2}}\int_0^t\tau^{1-\f{\t s_2(m+2)}{2}}\|f(\tau,\cdot)\|_{H^s}d\tau\\
&\le Ct^{2-\f{p_1(m+2)}{2}}\|f\|_{L^{\infty}([0, T],H^s)}.\tag3.29
\endalign
$$

Substituting (3.28)-(3.29) into (3.27) yields (3.25) for $0\le
p_1<p_1(m)$.

Next, we show (3.26).

Due to
$$
\p_tu^\wedge(t,\xi)=\int_0^t \bigg(\p_tV_2(t,|\xi|)V_1(\tau,|\xi|)
-\p_tV_1(t,|\xi|)V_2(\tau,|\xi|)\bigg)f^\wedge(\tau,\xi)d\tau,
$$

then one has  by  Minkowski inequality
$$\align \|\p_t&u(t,\cdot)\|_{H^{s-\f{m}{2(m+2)}+p_2}}\\
&\leq \int_0^t\bigg(\int_{\Bbb
R^n}\bigg|(1+|\xi|^2)^{\f{s}{2}-\f{m}{4(m+2)}+\f{p_2}{2}}
(\p_tV_2(t,|\xi|)V_1(\tau,|\xi|)
-\p_tV_1(t,|\xi|)V_2(\tau,|\xi|))f^\wedge(\tau,\xi)\bigg|^2d\xi\bigg)^{\f{1}{2}}d\tau\\
&\leq\int_0^t\|(1+|\xi|^2)^{\f{s}{2}-\f{m}{4(m+2)}+\f{p_2}{2}}
\p_tV_2(t,|\xi|)V_1(\tau,|\xi|)f^\wedge(\tau,\xi)\|_{L^2}d\tau\\
&\quad
+\int_0^t\|(1+|\xi|^2)^{\f{s}{2}-\f{m+4}{4(m+2)}+\f12\min\{\f{m+4}{2(m+2)},
\f{4}{m+2}\}-\f{p_2(m)-p_2}{2}}
\p_tV_1(t,|\xi|)V_2(\tau,|\xi|)f^\wedge(\tau,\xi)\|_{L^2}d\tau\\
\equiv& I_3+I_4.\tag3.30\endalign$$

Applying for Lemma 3.2  yields for $0<t\leq T$ and $0<\tau\leq T$
$$\align \quad I_3
\leq &C\int_0^t\|(1+|\xi|^2)^{\f{p_2}{2}}V_1(\tau,|\xi|)
(1+|\xi|^2)^{\f{s}{2}}f^\wedge(\tau,\xi)\|_{L^2}d\tau\\
\le&C\int_0^t\tau^{-\f{m+2}{2}p_2}\|f(\tau,\cdot)\|_{H^s}d\tau\\
\le&C_{p_2}t^{1-\f{m+2}{2}p_2}\|f\|_{L^{\infty}([0,T],H^s)}\tag3.31\endalign$$
and
$$
\align I_4&\leq C\int_0^t\|(1+|\xi|^2)^{\f12\min\{\f{m+4}{2(m+2)},
\f{4}{m+2}\}-\f{p_2(m)-p_2}{2}}V_2(\tau,|\xi|)
(1+|\xi|^2)^{\f{s}{2}}f^\wedge(\tau,\xi)\|_{L^2}d\tau\\
&\le
C\int_0^t\tau^{1-\f{m+2}{2}p_2}\|f(\tau,\cdot)\|_{H^s}d\tau\\
&\le
C_{p_2}t^{2-\f{m+2}{2}p_2}\|f\|_{L^{\infty}([0,T],H^s)}.\tag3.32
\endalign$$

Substituting (3.31)-(3.32) into (3.30) yields (3.26).

Consequently, Lemma 3.4 is proved. \hfill \qed

\vskip 0.3 true cm \centerline{\bf $\S 4$. Conormal spaces and
commutator relations} \vskip 0.3 true cm

In this section, we will give the definitions of  conormal spaces
related to our problems. To this end, as the first step, we look for
the basis of vector fields tangent to some surface (or surfaces).

{\bf Lemma 4.1.} {\it Let $\Gamma_0=\{(t,x): t\ge 0,
|x|^2=\ds\f{4t^{m+2}}{(m+2)^2}\}$,  then a basis of the $C^{\infty}$
vector fields tangent to $\Gamma_0$ is given by
$$\align &L_0=2t\p_t+(m+2)(x_1\p_1+\cdots+x_n\p_n);\\
&L_{i}=2t^{m+1}\p_i+(m+2)x_i\p_t,\quad i=1,2,\cdots,n;\\
&L_{ij}=x_i\p_j-x_j\p_i,\quad 1\le i<j\le n.\endalign$$}

{\bf Proof.} At first, we prove such an assertion:

{\it Given a smooth function $c(t,x)$ vanishing on $\Gamma_0$, then
there exists a smooth function $d(t,x)$ such that
$$c(t,x)=d(t,x)(|x|^2-\f{4t^{m+2}}{(m+2)^2}).\tag 4.1$$}

Indeed, it follows from Malgrange Preparation Theorem (see Theorem
7.5.6 of [16]) that there exist smooth functions $c_1(t, \t x)$,
$c_2(t,\t x)$ with $\t x=(x_1,..., x_{n-1})$ and $d(t,x)$ such that
$$c(t,x)=d(t,x)(|x|^2-\f{4t^{m+2}}{(m+2)^2})+x_nc_1(t,\t x)+c_2(t,\t x).$$

For $(t,x)\in\Gamma_0$, we have
$$0=\pm c_1(t,\t x)\sqrt{\f{4t^{m+2}}{(m+2)^2}-|\t x|^2}+c_2(t,\t x).$$

This yields $c_1(t,\t x)=c_2(t, \t x)=0$. Hence, we complete the
proof on (4.1).

Next we use the induction method and (4.1) to prove Lemma 4.1.

For $n=1$, we assume that the vector field
$L=a(t,x_1)\p_1+b(t,x_1)\p_{t}$ is tangent to
$\G_0^1\equiv\{x_1^2=\ds\f{4t^{m+2}}{(m+2)^2}\}$, which means
$$L(x_1^2-\f{4t^{m+2}}{(m+2)^2})=2a(t,x_1)x_1-\f{4b(t,x_1)}{m+2}t^{m+1}\equiv 0
\q\text{on $\G_0^1$}.\tag4.2$$

By (4.1), we know that there exists a smooth function $d(t,x_1)$
such that
$$2a(t,x_1)x_1-\f{4b(t,x_1)}{m+2}t^{m+1}=d(t,x_1)(x_1^2-\f{4t^{m+2}}{(m+2)^2})
.\tag4.3$$

This derives
$$b(t,0)=\f{d(t,0)t}{m+2}.\tag4.4$$

On the other hand, it follows  from (4.3) that there exist two
smooth functions $b_1(t,x_1)$ and $d_1(t,x_1)$ such that
$$2a(t,x_1)x_1-\f{4b(t,0)+4b_1(t,x_1)x_1}{m+2}t^{m+1}=(d(t,0)+d_1(t,x_1)x_1)\biggl(x_1^2-\f{4t^{m+2}}{(m+2)^2}\biggr)
.\tag4.5$$

Substituting (4.4) into (4.5) yields
$$a(t,x_1)=\f{2b_1(t,x_1)}{m+2}t^{m+1}+\ds\f{d(t,0)}{2}x_1+\ds\f{d_1}{2}\biggl(x_1^2-\f{4t^{m+2}}{(m+2)^2}\biggr)
.\tag4.6$$

Therefore,
$$\align L&=a(t,x_1)\p_1+(b(t,0)+b_1(t,x_1)x_1)\p_{t}\\
&=\ds\f{b_1}{m+2}L_1+\f{d(t,0)}{2(m+2)}L_0+\f{d_1}{2(m+2)^2}\biggl((m+2)x_1L_0-2tL_1\biggr)\\
&=\biggl(\ds\f{b_1}{m+2}-\f{d_1t}{(m+2)^2}\biggr)L_1+
\biggl(\ds\f{d_1x_1}{2(m+2)}+\f{d(t,0)}{2(m+2)}\biggr)L_0.
\endalign
$$

This yields the case of $n=1$ in Lemma 4.1.

By the induction hypothesis, we assume that Lemma 4.1  holds for
$n-1$.

We now show the case of $n$.

Assume that $L=a_0(t,x)\p_t+\ds\sum_{i=1}^na_i(t,x)\p_i$ is tangent
to $\Gamma_0$, and we  rewrite $L$ as
$$\align L &=\sum _{i=0}^na_i(t,\t x,0)\p_i
+\sum_{i=0}^nb_i(t,\t x,x_n)x_n\p_i\\
&\equiv M_{n-1}+a_n(t,\t x,0)\p_n +\sum_{i=0}^nb_i(t,\t
x,x_n)x_n\p_i,\endalign$$ where $\t x=(x_1, ..., x_{n-1})$,
$\p_0=\p_t$, and $M_{n-1}=\dsize\sum _{i=0}^{n-1}a_i(t, \t
x,0)\p_i$.

We can assert that $M_{n-1}$ is tangent to the surface $\{|\t
x|^2=\ds\f{4t^{m+2}}{(m+2)^2}\}$.

Indeed, due to
$$L(|x|^2-\f{4t^{m+2}}{(m+2)^2})=M_{n-1}(|\t x|^2-\f{4t^{m+2}}{(m+2)^2})+x_n\biggl(
2a_n(t,\t
x,0)+\ds\sum_{i=0}^nb_i\p_i(|x|^2-\f{4t^{m+2}}{(m+2)^2})\biggr),$$
and $L(|x|^2-\ds\f{4t^{m+2}}{(m+2)^2})\equiv 0$ on
$\{|x|^2=\ds\f{4t^{m+2}}{(m+2)^2}\}$, then
$L(|x|^2-\ds\f{4t^{m+2}}{(m+2)^2})\equiv 0$ holds true on $\{|\t
x|^2-\ds\f{4t^{m+2}}{(m+2)^2}=0, x_n=0\}$ and further $M_{n-1}(|\t
x|^2-\ds\f{4t^{m+2}}{(m+2)^2})\equiv 0$ on $\{|\t
x|^2=\ds\f{4t^{m+2}}{(m+2)^2}\}$ is derived.

By the induction hypothesis, we know that $M_{n-1}$ can be expressed
as a linear combination of
$L_0^{(n-1)}=2t\p_t+(m+2)(x_1\p_1+\cdots+x_{n-1}\p_{n-1})$,
$L_{i}^{(n-1)}=2t^{m+1}\p_i+(m+2)x_i\p_t$ ($i=1,2,\cdots,n-1$) and
$L_{ij}^{(n-1)}=x_i\p_j-x_j\p_i$ ($1\le i<j\le n-1$). On the other
hand, we have $x_n\p_t=\f{1}{m+2}L_0^{(n)}-\f{2}{m+2}t^{m+1}\p_n$
and $x_n\p_i=L_{ni}^{(n)}+x_i\p_n$. Thus, we can arrive at

$$L=\t a(t,x)\p_n+p(L_0^{(n)},L_i^{(n)}, L_{ij}^{(n)})_{1\le i<j\le n},\tag4.7$$
where $\t a(t,x)$ is a smooth function, $p(L_0^{(n)},L_i^{(n)},
L_{ij}^{(n)})_{1\le i<j\le n}$ represents a first order polynomial
of $L_0^{(n)}$, $L_i^{(n)}, ... ,L_{ij}^{(n)}$ with $1\le i\le n$
and $1\le i<j\le n$ respectively.

Since $L$ and $L_0^{(n)}, L_i^{(n)}, ... , L_{ij}^{(n)}$ are all
tangent to $\G_0$, then one has $\t a(t,x)\equiv 0$ on $\Gamma_0$.
This yields that there  exists a smooth function $d(t,x)$ such that
$\t a(t,x)=d(t,x)((m+2)^2|x|^2-4t^{m+2})$. It is noted that
$((m+2)^2|x|^2-4t^{m+2})\p_n=(m+2)x_nL_0^{(n)}-2tL_n^{(n)}
-(m+2)^2\dsize\sum_{i=1}^{n-1}x_iL_{ni}^{(n)}$. This, together with
(4.7), yields the proof on Lemma 4.1.

In order to apply for the commutator argument to treat our
degenerate equation whose characteristic cone and characteristic
surfaces have cusp singularities, we will use the following revised
vector fields tangent to $\G_0$:
$$\align &L_0=2t\p_t+(m+2)(x_1\p_1+\cdots+x_n\p_n);\\
&\bar L_{i}=2t^{\f{m}{2}+1}\p_i+(m+2)\f{x_i}{t^{\f{m}{2}}}\p_t,\quad i=1,2,\cdots,n;\\
&L_{ij}=x_i\p_j-x_j\p_i,\quad 1\le i<j\le n.\endalign$$

Let $[A,B]=AB-BA$ denote  the commutator of $A$ and $B$. By a direct
computation, one has

\vskip .2cm

 {\bf Lemma 4.2.} {\it For $1\le i\le n$ and $1\le i<j\le
n$,
$$\align &[L_0,\bar L_i]=0; \q [L_0,L_{ij}]=0;\\
&[\bar L_i,\bar L_j]=2(m+1)(m+2)L_{ij}+\f{m(m+2)}{2}(\f{x_j}{t^{\f{m}{2}+1}}\bar L_i-
\f{x_i}{t^{\f{m}{2}+1}}\bar L_j); \q [\bar L_i,L_{ij}]=\bar L_j; \\
&[\bar L_k,L_{ij}]=0\q \text{for $k\neq i$ and $k\neq j$};\\
&[L_{ij},L_{kl}]=0\q \text{for $1\le k<l\le n$, $k\neq i, l\neq j$};
\q [L_{ij},L_{ik}]=L_{kj}
\q\text{for $k\neq j$}.\\
\endalign
$$

In addition, let $P=\p_t^2-t^m\Delta$, then
$$
[P,L_0]=4P,  \q[P,L_{ij}]=0, \q [P,\bar L_{i}]=-m(m+2)
\f{x_i}{t^{\f{m}{2}+1}} P+ \f{m(m+2)}{4 t^2}\bar{L}_i.$$} {\bf
Remark 4.1.} {\it If we choose the smooth vector fields
$L_i=2t^{m+1}\p_i+(m+2)x_i\p_t$ instead of $\bar L_i$ in Lemma 4.2,
then a direct computation yields that for $1\le i\le n$
$$[P, L_i]=mt^{m-1}\biggl(\p_i L_0+(m+2)\ds\sum_{j\neq
i}\p_jL_{ij}+(n(m+2)-2)\p_i\biggr)$$ and
$$\text{$[P, L_i^m]$ will include the term $((m+2)x_i)^{m-1}m!\biggl(\p_i L_0+(m+2)\ds\sum_{j\neq
i}\p_jL_{ij}+(n(m+2)-2)\p_i\biggr)$}.$$

When the conormal
regularities of solution $u(t,x)$ to problem (1.1) are studied, we
will meet the following problem
$$PL_i^mu=L_i^m(f(t,x,u))+[P, L_i^m]u.$$
By Proposition 3.3, under the assumptions $(A_2)-(A_3)$, one can only
expect $u(t,x)\in C([0, T], H^{\f12+\f{2}{m+2}-})$, which implies
$PL_i^mu=L_i^m(f(t,x,u))+\text{\bf some terms of $C([0, T],
H^{\f{2}{m+2}-\f12-}(\Bbb R^n))$}$ and thus only $L_i^mu\in C([0,$
$T], H^{\f{6}{m+2}-\f12-}(\Bbb R^n))$ can be expected by Lemma 3.4.
Hence, for large $m$, one just only obtains $L_i^mu\in C([0, T],$
$H^{-\f12}(\Bbb R^n))$, which leads to the loss of regularities of
$L_i^mu$ and more losses of regularities of $\{L_i^lu\}_{l>m}$ can
be produced. It is noted that $L_i^ku\in C([0, T], L^2(\Bbb R^n))$
with $k\ge m$ should be obtained in order to show Theorem 1.1.
Hence, one can not use the smooth vector fields and commutator
arguments directly to show Theorem 1.1.}

Next, we introduce the vector fields tangent to
$\G_1^{\pm}=\{(t,x)\in \Bbb R_+\times\Bbb R^n:
x_1=\pm\ds\f{2t^{\f{m+2}{2}}}{m+2}\}$.

{\bf Lemma 4.3.} {\it The following vector fields are tangent to
$\Gamma_1^{\pm}$
$$\bar L_0=2t\p_t+(m+2)x_1\p_1,\quad \bar L_1=2t^{\f{m}{2}+1}\p_1+(m+2)\f{x_1}{t^{\f{m}{2}}}\p_t,
\quad R_k=\p_k, 2 \le k \le n.$$ Moreover, we have the following
commutator relations:
$$
\align &[\bar L_0, \bar L_1]=0;\quad [\bar L_0,
R_k]=-(m+2)R_k\quad\text{for
$k\ge 2$};\quad [\bar L_1, R_k]=0\quad\text{for $k\ge 2$};\\
&[P, \bar L_0]=4P+2mt^m\ds\sum_{i=2}^n\p_iR_i;\quad [P, \bar L_1]=-m(m+2)
\f{x_1}{t^{\f{m}{2}+1}} P+ \f{m(m+2)}{4 t^2}\bar{L}_1;\quad [P,
R_k]=0\quad\text{for $k\ge 2$}.
\endalign$$}

\vskip .1cm

{\bf Proof.} This can be verified directly, we omit it here. \qquad\qquad\qquad\qquad
\qquad\qquad\qquad\qquad\qquad\qquad\qquad\hfill
\qed

\vskip .2cm

{\bf Remark 4.2.} {\it In the expression of commutator $[P, \bar
L_1]$, there appear a singular factor $\ds\f{1}{t^2}$ before $\bar
L_1$. This will produce such an equation on $\bar L_1u$ from (1.1):
$$P\bar L_1u-\f{m(m+2)}{4t^2}\bar L_1u-f_u'(t,x,u)\bar L_1u=
(\bar L_1f)(t,x,u)-m(m+2) \f{x_1}{t^{\f{m}{2}+1}}f(t,x,u).$$ Such a
degenerate equation with a singular coefficient $\ds\f{1}{t^2}$ has
a bad behavior near $t=0$, and thus it is not suitable to use the
commutator argument on $\bar L_1u$ (or more general $\bar L_1^ku$)
to derive the regularity of $\bar L_1u$(or $\bar L_1^ku$).}

Based on the preparations above, we will define the conormal spaces
which will be required later on. To this end, such terminologies as
in [2]-[3] are introduced:

$\{M_1, \cdot\cdot\cdot, M_k\}$ stands for a collection of vector
fields with bounded coefficients on an open set $\O\subset\Bbb R^n$
such that all commutators $[M_i, M_j]$ are in the linear span over
$C^{\infty}(\O)$ of $M_1, \cdot\cdot\cdot, M_k$.

\vskip .2cm

{\bf Definition 4.1} ({\bf Admissible function}) {\it A function
$h(x)\in L^{\infty}(\O)\cap C^{\infty}(\O)$  is called admissible
with respect to $\{M_1, ..., M_k\}$ if $M_1^{j_1}\cdot\cdot\cdot
M_k^{j_k}h\in L^{\infty}(\O)\cap C^{\infty}(\O)$ for all $(j_1, ...,
j_k)$.}

Obviously, the linear span of $\{M_1, ..., M_k\}$ with admissible
coefficients is a Lie algebra of vector fields on $\O$.

We now define the admissible tangent  vector  fields related to the
surface $\G_0$.

\vskip .2cm

{\bf Definition 4.2} ({\bf Admissible tangent vector  fields of
$\G_0$}) {\it

(1) Let $\O_1$ be a region of the form $\{(t,x):0<t<C|x|\le \ve\}$
and $\Cal{S}_1$ be the Lie algebra of vector fields with admissible
coefficients on $\O_1$ generated by $\{|x|\p_t, t^{\f{m}{2}}|x|\p_i,
L_{ij}, i,j=1,2,\cdots,n\}$.

(2) Let $\O_2$ be a region of the form
$\{(t,x):|x|<Ct\le\ve\}\cap\{(t,x):||x|-\f{2}{m+2}t^{\f{m+2}{2}}|<Ct^{\f{m+2}{2}}\}$
and $\Cal{S}_2$ be the Lie algebra of vector fields with admissible
coefficients on $\O_2$ generated by $\{L_0, \bar L_i,
L_{ij},i,j=1,2,\cdots,n\}$ in Lemma 4.2.

(3) Let $\O_3$ be a region of the form
$\{(t,x):|x|<Ct\le\ve\}\cap\{(t,x):
t^{\f{m+2}{2}}<C||x|\pm\f{2}{m+2}t^{\f{m+2}{2}}|\}$ and $\Cal{S}_3$
be the Lie algebra of vector fields with admissible coefficients on
$\O_3$ which generated by $\{t\p_t, t^{m+1}\p_i, L_{ij}, i,j$ $=1,2,
\cdots, n\}$.}

Next, the conormal space $I^{\infty}H^{s}(\G_0)$ with $0\le
s<\f{n}{2}$ is defined.

\vskip .2cm

{\bf Definition 4.3} ({\bf Conormal space $I^{\infty}H^{s}(\G_0)$}).
{\it Define the function $u(t,x)\in I^{\infty}H^{s}(\G_0)$ in
$\{(t,x): 0\le t\le T,$ $x\in\Bbb R^n\}$ if, away from
$\{|x|=t=0\}$, $Z_1\cdots Z_ku\in L^{\infty}([0, T],
H^{s}(\Bbb R^{n}))$ for any $k\in\Bbb N\cup\{0\}$ and all
smooth vector fields $Z_1,\cdots,Z_j\in \{L_0, \bar L_i,
L_{ij},i,j=1,2,\cdots,n\}$, and near $\{|x|=t=0\}$, the following
properties hold:

(1) If $h_1(t,x)\in C^{\infty}(\Bbb R^{n+1}\setminus\{0\})$ is
homogeneous of degree zero and supported on $\O_1=\{(t,x): 0\le
t<C|x|\le\ve\}$, then $Z_1\cdots Z_k(h_1(t,x)u(t,x))\in
L^{\infty}([0, T], H^{s}(\Bbb R^{n}))$ for all $Z_1,\cdots, Z_k\in
\Cal{S}_1$.

(2) If $h_2(t,x)\in C^{\infty}(\Bbb R^{n+1}\setminus\{0\})$ is
homogeneous of degree zero and supported on $\{(t,x):|x|<Ct\le\ve\}$
and $\chi(\theta)\in C^{\infty}$ has compact support near
$\{\theta=1\}$, then $Z_1\cdots
Z_k(h_2(t,x)\chi(\f{(m+2)|x|}{2t^{\f{m+2}{2}}})u)\in L^{\infty}([0,
T], H^{s}(\Bbb R^{n}))$ for all $Z_1,\cdots, Z_k\in \Cal{S}_2$.

(3) If $h_3(t,x)\in C^{\infty}(\Bbb R^{n+1}\setminus\{0\})$ is
homogeneous of degree zero and supported on $\{(t,x):|x|<Ct\le\ve\}$
and $\chi_0(\theta)\in C^{\infty}$ has compact support away
$\{\theta=1\}$, then $Z_1\cdots
Z_k(h_3(t,x)\chi_0(\f{(m+2)|x|}{2t^{\f{m+2}{2}}})u)\in
L^{\infty}([0, T], H^{s}(\Bbb R^{n}))$ for all $Z_1,\cdots, Z_j\in
\Cal{S}_3$.}

It is noted that  $h_1(t,x)$,
$h_2(t,x)\chi(\f{(m+2)|x|}{2t^{\f{m+2}{2}}})$, and
$h_3(t,x)\chi_0(\f{(m+2)|x|}{2t^{\f{m+2}{2}}})$ are admissible
functions on domains $\O_1$, $\O_2$ and $\O_3$ respectively,
moreover are in the space $L^{\infty}([0, \infty),
H^{\f{n}{2}-}(\Bbb R^n))$.

Because some vector fields (for examples, $\bar L_i$, $i=1, ..., n$)
in Definition 4.3 has no good commutator relations (i.e., the
coefficients of commutator are not admissible, one can also see the
explanations in Remark 4.2) with  $P=\p_t^2-t^m\Delta$,
we have to look for some auxilliary relations among those vector fields
which possess good commutator relations (for examples, $L_0, L_{ij}$
in Lemma 4.2, and $\bar L_0, R_j (2\le j\le n)$ in Lemma 4.3) and
$P$.

Formally, it follows from a direct computation that

$$\cases &\p_i=\ds\f{\Big(4t^{m+2}-(m+2)^2\ds\sum_{j\neq i}x_j^2 \Big)\bar L_i
+(m+2)^2\sum_{j\neq i}x_ix_j\bar L_j-2(m+2)x_it^{\f{m}{2}+1}L_0}
{2t^{\f{m}{2}+1}(4t^{m+2}-(m+2)^2|x|^2)},\quad i=1, ..., n;\\
&\p_{t}=\ds\f{2 t^{m+1}L_0-(m+2)\ds\sum_{i=1}^nx_it^{\f{m}{2}}\bar
L_i}{4t^{m+2}-(m+2)^2|x|^2}.\endcases\tag4.8$$

According to this and some crucial observations, we have

{\bf Lemma 4.4.} {\it Let $\O_i (i=1,2,3)$ be given in Definition
4.2, one has

(1) On $\O_1$, set $N_1^0=|x|\p_{t}$, $N_1^i=t^{\f{m}{2}}|x|\p_i$
with $i=1, ..., n$, then
$$\align (N_1^0)^2=&\f{1}{4t^{m+2}-(m+2)^2|x|^2}
\bigg(-4|x|^2t^{m+2}P-|x|^2t^m \sum_{j=1}^n\bar L_j^2+4|x|t^{m+1}N_1^0L_0\\
&\hskip 3cm  +(m+1)|x|^2t^mL_0 +\Big((2nm+2n-2m-2)t^{m+1}|x|-\f{m(m+2)^2 |x|^3}{2t}\Big)N_1^0\bigg);\\
(N_1^i)^2=&\f{1}{4t^{m+2}-(m+2)^2|x|^2}
\bigg(-|x|^2(4t^{m+2}-(m+2)^2\sum_{j\neq i}x_j^2))P+t^m|x|^2L_0^2
-|x|^2t^m\sum_{j\neq
i}\bar L_j^2\\
& \quad -2(m+2)x_i|x|t^{\f{m}{2}} N_1^iL_0+\Big( ((n-1)(m+2)-2)
t^{m+2}
-\f{m(m+2)^2 }{4}\sum_{j \ne i}x_j^2\Big) \f{|x|^2}{t^2}L_0 \bigg)\\
& \quad + a_1^i N_1^i + \sum_{j \ne i} b_1^j N_1^j,\endalign$$ where
$a_1^i, b_1^j$ are admissible on $\O_1$.

 (2) On $\O_2$,  set $N_2^i=(|x|-\f{2}{m+2}t^{\f{m+2}{2}})\p_i$
with $i=1, ..., n$, then
$$
\align \bar
L_i&=\f{1}{2t^{\f{m}{2}+1}}\biggl((m+2)x_iL_0-(m+2)^2\ds\sum_{k\not=i}x_kL_{ik}\\
&\qquad -(m+2)((m+2)|x|+2t^{\f{m}{2}+1})N_2^i\biggr)\tag4.9
\endalign
$$
and
$$\align (N_2^i)^2=&\f{1}{(m+2)((m+2)|x|+2t^{\f{m+2}{2}})}
\bigg(-(|x|-\f{2}{m+2}t^{\f{m+2}{2}})(4t^2-(m+2)^2\f{\ds\sum_{j\neq
i}x_j^2}{t^m})P\\
&+(|x|-\f{2}{m+2}t^{\f{m+2}{2}})L_0^2
-(|x|-\f{2}{m+2}t^{\f{m+2}{2}})\sum_{j\neq
i}\bar L_j^2\\
&-2(m+2)x_i N_2^iL_0
+(|x|-\f{2}{m+2}t^{\f{m+2}{2}})\Big( (n-1)(m+2)-2-\f{m(m+2)^2}{4t^{m+2}} \sum_{j \ne i} x_j^2\Big)L_0\bigg)\\
&+ a_2^i N_2^i + \sum_{j \ne i} b_2^j N_2^j,\tag4.10\endalign$$
where $a_2^i$ and $b_2^j$ are admissible on $\O_2.$

 Thus, one has from (4.9) and
(4.10) that for $i=1, ..., n$,
$$
\align (N_2^i)^2=&a_0P+a_1L_0^2+\ds\sum_{\Sb 1\le i<k\le
n,\\
1\le m<l\le n\endSb}a_{ik}^{ml}L_{ik}L_{ml} +\ds\sum_{1\le i<k\le
n}b_{ik}L_0L_{ik}+\ds\sum_{1\le i\le n}b_iN_2^iL_0\\
&+\ds\sum_{\Sb 1\le i\le n,\\ 1\le m<l\le
n\endSb}b_i^{ml}N_2^iL_{ml}+\ds\sum_{1\le i\le
n}c_{0i}N_2^i+cL_0,\tag4.11
\endalign$$
where the coefficients $a_0, a_1, a_{ik}^{ml}, b_{ik}, b_i,
b_i^{ml}, c_{0i}$ and $c$ are admissible on $\O_2$.

(3) On $\O_3$, set $N_3^0=t\p_t$, $N_3^i=t^{\f{m}{2}+1}\p_i$ with
$i=1, ..., n$, then

$$\align (N_3^0)^2=&\f{1}{4t^{m+2}-(m+2)^2|x|^2}\bigg(-4t^{m+4}P
-t^{m+2}\sum_{j=1}^n\bar L_j^2+4t^{m+2}N_3^0L_0\\
&+(m+2) t^{m+2}L_0 +\Big( 2(n-1)(m+2) t^{m+2} -\f{(4+m)(m+2)^2
|x|^2}{2}\Big)N_3^0\bigg)\endalign$$
and
$$\align (N_3^i)^2=&\f{1}{(m+2)^2|x|^2-4t^{m+2}}
\bigg(t^{m+2}(4t^{m+2}-(m+2)^2\sum_{j\neq
i}x_j^2)P-t^{2(m+1)}L_0^2+t^{2m+2}\sum_{j\neq i}\bar L_j^2\\
& \quad +2(m+2)x_it^{\f{3m+2}{2}}N_3^iL_0 \Big(
(2-(m+2)(n-1))t^{2m+2} +
\f{m(m+2)^2}{4}t^m\sum_{j \ne i}x_j^2\Big) L_0\bigg) \\
&+ a_3^i N_3^i + \sum_{j \ne} b_3^j N_3^j,\endalign$$} where $a_3^i,
b_3^j$ are admissible on $\O_3$.

 {\bf Remark 4.3.} {\it It is also easy to verify that
the coefficients on each region $\O_i$ $(i=1,3)$ in Lemma 4.4  are
admissible with respect to the vector fields $\Cal{S}_i$
respectively.}

We start to define the admissible tangent vector fields related to
the surfaces $\G_1^+$ and $\G_1^-$.

\vskip .2cm

{\bf Definition 4.4} ({\bf Admissible tangent vector  fields of
$\G_1^{\pm}$}) {\it

(1) Let $W_1$ be a region of the form $\{(t,x): 0<t<C|x_1|\le\ve\}$
and $\Cal{M}_1$ be the Lie algebra of vector fields with admissible
coefficients on $W_1$ generated by $\{x_1\p_1, x_1\p_t, \p_i,
i=2,\cdots,n\}$.

(2) Let $W_{2,\pm}$ be a region of the form
$\{(t,x):|x_1|<Ct<\ve\}\cap\{(t,x):|x_1\mp\f{2}{m+2}t^{\f{m+2}{2}}|<Ct^{\f{m+2}{2}}\}$
and $\Cal{M}_{2\pm}$ be the Lie algebra of vector fields with
admissible coefficients on $W_{2,\pm}$ generated by $\{\bar L_0,
\bar L_1, R_2, \cdots, R_n\}$ in Lemma 4.3.

(3) Let $W_3$ be a region of the form
$\{(t,x):|x_1|<Ct\le\ve\}\cap\{(t,x):
t^{\f{m+2}{2}}<C|x_1\pm\f{2}{m+2}t^{\f{m+2}{2}}|\}$ and $\Cal{M}_3$
be the Lie algebra of vector fields with admissible coefficients on
$W_3$ generated by $\{t\p_t, t^{\f{m+2}{2}}\p_1, \p_i, i=2, \cdots,
n\}$.}

{\bf Remark 4.4.} {\it On $W_{2,\pm}$, for the convenience of
computation, we sometimes use the equivalent vector fields
$\Cal{M}_{2\pm}=\{\bar L_0, N_{2,\pm}, R_2, \cdots, R_n\}$ with
$N_{2,\pm}\equiv (x_1\mp\f{2}{m+2}t^{\f{m+2}{2}})\p_1$ instead of
$\{\bar L_0, \bar L_1, R_2, \cdots, R_n\}$ from now on. The
equivalence comes from the following facts:
$$\align
&N_{2,\pm}=\f{t^{\f{m}{2}+1}}{(m+2)^2(x_1\pm\f{2}{m+2}t^{\f{m+2}{2}})}
\biggl(\f{(m+2)x_1}{t^{\f{m}{2}+1}}\bar L_0-2\bar L_1\biggr),\\
&\bar
L_0=\f{t^{\f{m}{2}+1}}{(m+2)x_1}\biggl(\f{(m+2)^2(x_1\pm\f{2}{m+2}t^{\f{m+2}{2}})}
{t^{\f{m}{2}+1}}N_{2,\pm}+2\bar L_1\biggr),
\endalign
$$
where all related coefficients are admissible on $W_{2,\pm}$.}

Similarly, we  define the conormal space
$I^{\infty}H^s(\G_1^+\cup\G_1^-)$ with $0\le s<\f{n}{2}$.

\vskip .2cm

{\bf Definition 4.5} ({\bf Conormal space}
$I^{\infty}H^s(\G_1^+\cup\G_1^-)$). {\it Define the function
$u(t,x)\in I^{\infty}H^s(\G_1^+\cup\G_1^-)$ in $t\ge 0$ if, away
from $\{x_1=t=0\}$, $Z_1\cdots Z_ku\in L^{\infty}([0, T], H^s(\Bbb
R^n))$ for all smooth vector fields $Z_1,\cdots,Z_k\in \{\bar L_0,
\bar L_1, R_2, \cdots, R_n\}$ in Lemma 4.3, and near $\{x_1=t=0\}$,
the following properties hold:

1) If $h_1(t,x_1)\in C^{\infty}(\Bbb R^{2}\setminus\{0\})$ is
homogeneous of degree zero and supported on $W_1=\{(t,x): 0\le
t<C|x_1|\le\ve\}$, then $Z_1\cdots Z_k(h_1(t,x_1)u(t,x))\in
L^{\infty}([0, T], H^s(\Bbb R^n))$ for all $Z_1,\cdots, Z_k\in
\Cal{M}_1$.

2) If $h_2(t,x_1)\in C^{\infty}(\Bbb R^{2}\setminus\{0\})$ is
homogeneous of degree zero and supported on
$\{(t,x_1):|x_1|<Ct\le\ve\}$ and $\chi_{\pm}(\theta)\in C^{\infty}$
has compact support near $\{\theta=\pm1\}$, then $Z_1\cdots
Z_k(h_2(t,x_1)\chi_{\pm}(\f{(m+2)x_1}{2t^{\f{m+2}{2}}})u)\in
L^{\infty}([0, T], H^s(\Bbb R^n))$ for all $Z_1,\cdots, Z_k\in
\Cal{M}_{2,\pm}$.

3) If $h_3(t,x_1)\in C^{\infty}(\Bbb R^{2}\setminus\{0\})$ is
homogeneous of degree zero and supported on
$\{(t,x):|x_1|<Ct\le\ve\}$ and $\chi_0(\theta)\in C^{\infty}$ has
compact support away $\{\theta=\pm1\}$, then $Z_1\cdots
Z_k(h_3(t,x_1)\chi_0(\f{(m+2)x_1}{2t^{\f{m+2}{2}}})u)\in
L^{\infty}([0, T], H^s(\Bbb R^n))$ for all $Z_1,\cdots, Z_k\in
\Cal{M}_3$.}

Obviously, the cutoff functions $h_1(t,x_1)$,
$h_2(t,x_1)\chi_{\pm}(\f{(m+2)x_1}{2t^{\f{m+2}{2}}})$, and
$h_3(t,x_1)\chi_0(\f{(m+2)x_1}{2t^{\f{m+2}{2}}})$ are admissible on
domains $W_1$, $W_{2,\pm}$ and $W_3$ respectively, moreover are in
the space $L^{\infty}([0, \infty), H^{\f{n}{2}-}(\Bbb R^n))$.

Similar to Lemma 4.4 and by some crucial observations, we have

\vskip .2cm

{\bf Lemma 4.5.}  {\it Let $W_1, W_{2,\pm}$ and $W_3$ be given in
Definition 4.4, one has

(1) On $W_1$, set $N_1=x_1\p_{t}$, then
$$\align N_1^2=&\f{1}{(m+2)^2x_1^2-4t^{m+2}}
\bigg((m+2)^2x_1^4P+x_1^2t^m\bar L_0^2-4x_1t^{m+1}N_1\bar L_0\\
&+ (m+2)^2 x_1^4t^m \ds\sum_{i=2}^n R_i^2-(m+2)x_1^2t^m\bar L_0
+2(m+4)x_1t^{m+1}N_1\bigg).\endalign$$

(2) On $W_{2,\pm}$, set $N_{2,\pm}=
(x_1\mp\f{2}{m+2}t^{\f{m+2}{2}})\p_1$, then
$$\align N_{2,\pm}^2=&\f{x_1\mp\f{2}{m+2}
t^{\f{m+2}{2}}}{(m+2)^2(x_1\pm\f{2}{m+2}t^{\f{m+2}{2}})}
\bigg(4t^2P-\bar L_0^2
+4t^{m+2}\ds\sum_{i=2}^nR_i^2+2\bar L_0\bigg)\\
&\qquad
+\f{2x_1}{(m+2)(x_1\pm\f{2}{m+2}t^{\f{m+2}{2}})}N_{2,\pm}\bar L_0
-\f{2(x_1 \pm
t^{\f{m+2}{2}})}{(m+2)(x_1\pm\f{2}{m+2}t^{\f{m+2}{2}})}N_{2,\pm}.\endalign$$

(3) On $W_3$, set $N_3=t\p_{t}$, $N_{3'}=t^{\f{m+2}{2}}\p_1$, then
$$\align N_3^2=&\f{1}{(m+2)^2x_1^2-4t^{m+2}}\bigg((m+2)^2x_1^2t^2P
+t^{m+2}\bar L_0^2-4t^{m+2}N_3\bar L_0\\
&+(m+2)^2x_1^2t^{m+2}\ds\sum_{i=2}^nR_i^2-(m+2)t^{m+2}\bar L_0
+\big((m+2)^2x_1^2+2(m+2)t^{m+2}\big)N_3\bigg)\endalign$$ and
$$\align N_{3'}^2=&\f{1}{(m+2)^2x_1^2-4t^{m+2}}
\bigg(4t^{m+4}P-t^{m+2}\bar L_0^2+2(m+2)x_1t^{\f{m+2}{2}}N_{3'}\bar L_0\\
&+4t^{2(m+2)}\ds\sum_{i=2}^nR_i^2+2t^{m+2}\bar L_0 -3(m+2)^2
x_1t^{\f{m+2}{2}}N_{3'}\bigg).\endalign$$}

{\bf Remark
4.5.} {\it As in Remark 4.3, one can easily  verify that the
coefficients on each domain in Lemma 4.5  are admissible with
respect to the corresponding vector fields.}

\vskip .2cm

 Finally, we define the conormal space
$I^kL^{\infty}_{loc}(\G_0\cup\G_1^{\pm}\cup\G_2^{\pm})$ of order
$k$, which are related to the surfaces $\G_0$, $\G_1^{\pm}$ and
$\G_2^{\pm}$ in $\Bbb R_+\times\Bbb R^2$. For this end, at first we
will introduce the admissible vector fields as in Definition 4.2 and
Definition 4.4.

Set
$$l_i^{\pm}=\G_i^{\pm}\cap \G_0\quad \text{for $i=1,2$},
\quad l_3^{\pm,\pm}=\G_1^{\pm}\cap \G_2^{\pm}.$$

For small fixed constant $\dl>0$ we define the following domains:
$$
\align
&\O_1^{\pm}=\{(t,x): t>0, |x_1\mp\ds\f{2t^{\f{m+2}{2}}}{m+2}|< \dl t^{\f{m+2}{2}}, 
\quad |x_2|<\dl t^{\f{m+2}{2}}\},\\
&\O_2^{\pm}=\{(t,x): t>0, |x_2\mp\ds\f{2t^{\f{m+2}{2}}}{m+2}|<\dl t^{\f{m+2}{2}}, 
\quad |x_1|<\dl t^{\f{m+2}{2}}\},\\
&\O_3^{\pm,\pm}\{(t,x): t>0, |x_1\mp\ds\f{2t^{\f{m+2}{2}}}{m+2}|<
\dl t^{\f{m+2}{2}},\quad
|x_2\mp\ds\f{2t^{\f{m+2}{2}}}{m+2}|< \dl t^{\f{m+2}{2}}\},\\
&\O_{i+3}^{\pm}=\{(t,x): t>0, |x_i\mp\ds\f{2t^{\f{m+2}{2}}}{m+2}|<
t^{\f{m+2}{2}}, \quad\text{away from the lines
$l_i^{\pm}$ and $l_3^{\pm,\pm}$}\},\quad i=1,2.\\
\endalign
$$
In addition,
$\O_1^{\pm}\cup\O_2^{\pm}\cup\O_3^{\pm,\pm}\cup\O_4^{\pm}\cup\O_5^{\pm}\cup(\cup_{j=1}^N\O_j)$
is an open cusp conic covering of ${\bar\Bbb R}_+^3\setminus\{O\}$
such that $\O_j$ ($1\le j\le N$) intersects at most one surface in
$\{\G_0, \G_1^{\pm}, \G_2^{\pm}\}$.

In $\O_1^{\pm}$, set ${\Cal M}_1^{\pm}=\{L_0, \bar L_1,
M_1^{\pm}\}$, here $ M_{1}^{\pm}=\ds\f{x_2}{t^{\f{m}{2}}}(\p_t\pm
t^{\f{m}{2}}\p_1)+( \f{2t^{\f{m}{2}+1}}{m+2}\mp x_1)\p_2$;

In $\O_2^{\pm}$, set ${\Cal M}_2^{\pm}=\{L_0, \bar L_2,
M_2^{\pm}\}$, here $M_{2}^{\pm}=\ds\f{x_1}{t^{\f{m}{2}}}(\p_t\pm
t^{\f{m}{2}}\p_2)+( \f{2t^{\f{m}{2}+1}}{m+2}\mp x_2)\p_1;$

In $\O_3^{\pm,\pm}$, set ${\Cal M}_3^{\pm,\pm}=\{L_0,N_1^{\pm, \pm},
N_2^{\pm, \pm}\}$, here $N_1^{\pm, \pm}=\ds\f{(\pm x_1\mp
x_2)}{t^{\f{m}{2}}}\p_t\pm (\ds\f{2t^{\f{m}{2}+1}}{m+2} \mp
x_2)\p_1\pm(\pm x_1-\ds\f{2t^{\f{m}{2}+1}}{m+2})\p_2$,
$N_2^{\pm,\pm}=t\p_t+t^{\f{m}{2}+1}(\pm\p_1\pm \p_2);$

In $\O_{i+3}^{\pm}$ (i=1,2), set ${\Cal M}_{i+3}^{\pm}=\{L_0,
R_{1,i}^{\pm}, R_{2,i}\}$, here $R_{1,i}^{\pm}=(x_i\mp
\ds\f{2t^{\f{m}{2}+1}}{m+2})\p_i,$ $R_{2, 1}=t\p_2$ or $R_{2,
2}=t\p_1;$

In $\O_{j}$ ($j=1, \cdot\cdot\cdot, N$), set ${\Cal M}_0=\{L_0, \bar
L_1, \bar L_2, L_{12}\}$.

Those vector fields with admissible coefficients generated by $\Cal
M_1^{\pm}, \Cal M_2^{\pm}, ..., {\Cal M}_0$ respectively are called
the admissible vector fields of surface variety $\{\G_0, \G_1^{\pm},
\G_2^{\pm}\}$.

Next, we give the definition of conormal space
$I^kL^{\infty}_{loc}(\G_0\cup\G_1^{\pm}\cup\G_2^{\pm})$.

\vskip .2cm

{\bf Definition 4.6} {\bf (Conormal space
$I^kL^{\infty}_{loc}(\G_0\cup\G_1^{\pm}\cup\G_2^{\pm})$)} {\it We
call a function $u(t,x)\in
I^kL^{\infty}_{loc}(\G_0\cup\G_1^{\pm}\cup\G_2^{\pm})$ if
$Z^{\al}\biggl(\chi(\ds\f{(m+2)x}{2t^{\f{m+2}{2}}})$
$u(t,x)\biggr)\in L^{\infty}_{loc}([0,T]\times\Bbb R^2)$ holds for
any $|\al|\le k$ and the homogeneous cut-off function
$\chi(\ds\f{(m+2)x}{2t^{\f{m+2}{2}}})$ of degree zero whose support
lies in some fixed conic neighborhood  of $\O_i^+ (i=1,2),
\O_3^{\pm,\pm},$ $\O_{i+3}^+ (i=1,2),\O_j (j=1,..., N)$, and $Z$
represents the admissible tangent vector in related domains.}

\vskip 0.3 true cm \centerline{\bf $\S 5$. Local existence of
solution to problem (1.1)} \vskip 0.3 true cm

In this section, we will show the local existence of the low
regularity solution to (1.1) in Theorem 1.1. At first, we study the
2-D case under the condition $(A_3)$ since the case of condition
$(A_1)$ is completely analogous and even simpler.

\vskip .2cm

{\bf Theorem 5.1.} {\it  Under the assumptions $(A_3)$ and $m\le 9$,
there exists a constant $T>0$ such that (1.1) has a local solution
$u\in L^{\infty}{([0, T]\times\Bbb R^2)}\cap C([0,T],
H^{\f{m+6}{2(m+2)}-}(\Bbb R^2))\cap C((0, T],$
$H^{\f{m+3}{m+2}-}(\Bbb R^2))\cap C^1([0,T], H^{\f{1}{m+2}-}(\Bbb
R^2))$.}

\vskip .1cm

{\bf Proof.}  Set
$$
\cases &\p_t^2u_1-t^{m}\triangle u_1=0,\\
&u_1(0,x)=0,\qquad \p_t{u_1}(0,x)=\vp(x).\endcases\tag5.1$$

Due to  $\vp(x)\in H^{\f12-}(\Bbb R^2)$ by Lemma 2.1.(ii), then it
follows from Lemma 2.4 and Proposition 3.3 that for any fixed
$0<\dl<\f{1}{2(m+2)}$
$$u_1(t,x)\in L^{\infty}([0, 1]\times\Bbb R^2)\cap C([0,1],
H^{\f{m+6}{2(m+2)}-\dl}(\Bbb R^2))\cap C((0, 1],
H^{\f{m+3}{m+2}-\dl}(\Bbb R^2))\cap C^1([0,1],
H^{\f{1}{m+2}-\dl}(\Bbb R^2)),$$ which satisfies for $t\in [0, 1]$
$$\|u_1(t,\cdot)\|_{L^{\infty}(\Bbb R^2)}+\|u_1(t,\cdot)\|_{H^{\f{m+6}{2(m+2)}-\dl}(\Bbb R^2)}+
t^{\f{m}{4}}\|u_1(t,\cdot)\|_{H^{\f{m+3}{m+2}-\dl}(\Bbb
R^2)}+\|\p_tu_1(t,\cdot)\|_{H^{\f{1}{m+2}-\dl}(\Bbb R^2)}\le
C(\dl).\tag5.2$$

Let
$$
\cases &\p_t^2u_2-t^m\Delta u_2=f(t,x,0),\\
&u_2(0,x)=\p_t{u_2}(0,x)=0.
\endcases\tag5.3
$$

Since $f(t,x,0)\in C^{\infty}([0, \infty)\times\Bbb R^2)$ has a
compact support on $x$, then (5.3) has a $C^{\infty}([0,
\infty)\times\Bbb R^2)$ solution $u_2(t,x)$ which possesses  a
compact support with respect to the variable $x$ when $t\ge 0$ is
fixed.

Set $v(t,x)=u(t,x)-u_1(t,x)-u_2(t,x)$, then it follows from (1.1),
(5.1) and (5.3) that
$$
\cases &\p_t^2v-t^{m}\triangle v=f(t,x,u_1+u_2+v)-f(t,x,0),\\
&v(0,x)=\p_t v(0,x)=0.\endcases\tag5.4$$

For  $w(t,x)\in C([0,T], H^{\f{m+6}{2(m+2)}+p_0(m)-\dl})\cap
C((0,T], H^{\f{m+6}{2(m+2)}+p_1(m)-\dl}) \cap C^1([0,T],$
$H^{\f{3}{m+2}+p_2(m)-\dl})$ satisfying for $t\in [0, T]$
$$
\align |\|w(t,\cdot)\||\equiv&\|w(t,\cdot)\|_{H^{\f{m+6}{2(m+2)}+
p_0(m)-\dl}(\Bbb R^2)}
+t^{\f{(m+2)p_1(m)}{2}-2}\|w(t,\cdot)\|_{H^{\f{m+6}{2(m+2)}+p_1(m)-\dl}(\Bbb
R^2)}\\
&+\|\p_tw(t,\cdot)\|_{H^{\f{3}{m+2}+
p_2(m)-\dl}(\Bbb R^2)}<\infty, \tag5.5\endalign$$ where $p_0(m)=\cases \ds\f{4}{m+2}
\quad\text{for $m\ge 2$}\\
1\quad\text{for $m\le 2$}\endcases$, $p_1(m)$ and $p_2(m)$ have been
defined in Lemma 3.4, and $0<T\le 1$, we define the set $G$ as
follows
$$
\align G\equiv \bigg\{w\in C([0,T], H^{\f{m+6}{2(m+2)}+
p_0(m)-\dl})\cap C((0,T], H^{\f{m+6}{2(m+2)}+p_1(m)-\dl}) &\cap
C^1([0,T],
H^{\f{3}{m+2}+p_2(m)-\dl}):\\
&\qquad \dsize\sup_{t\in [0,T]}|\|w(t,\cdot)\||\leq 1\bigg\}.
\endalign$$

Denote by
$$
\align &E(f(t,x,u)-f(t,x,0))\equiv\\
&\quad \biggl(\int_0^t
(V_2(t,|\xi|)V_1(\tau,|\xi|)-V_1(t,|\xi|)V_2(\tau,|\xi|))
(f(\tau,x,u(\tau,x))-f(\tau,x,0))^\wedge(\xi)d\tau\biggr)^\vee(t,x)
\endalign$$
and define a nonlinear mapping $\Cal F$ as follows
$$\Cal F(w)=E(f(t,x,u_1+u_2+w)-f(t,x,0)).\tag5.6$$

We now show that the mapping $\Cal F$ is from $G$ into itself and is
contractible for small $T$.

By (3.25) in Lemma 3.4 (taking $s_0=p_0(m)-\f{\dl}{2}$)  and (5.2),
we have for $w\in G$
$$
\align \|\Cal F(w)(t,\cdot)\|_{H^{\f{m+6}{2(m+2)}+p_0(m)-\dl}}&\le
Ct^{2-\f{(m+2)p_0(m)}{2}+\f{(m+2)\dl}{4}}
\|f(t,\cdot,u_1(t,\cdot)+u_2(t,\cdot)+w(t,\cdot))-f(t,\cdot,0)\|_{H^{\f{m+6}{2(m+2)}
-\f{\dl}{2}}}\\
&\le
Ct^{2-\f{(m+2)p_0(m)}{2}+\f{(m+2)\dl}{4}}\|u_1(t,\cdot)+u_2(t,\cdot)
+w(t,\cdot)\|_{H^{\f{m+6}{2(m+2)}-\f{\dl}{2}}}\\
&\le C(\dl)t^{2-\f{(m+2)p_0(m)}{2}+\f{(m+2)\dl}{4}},\tag5.7
\endalign$$
here we have used the follows facts:

$f(u)\in L^{\infty}([0, T]\times\Bbb R^n)\cap L^{\infty}([0,
T],H^s(\Bbb R^n))$ if $u\in L^{\infty}([0, T]\times\Bbb R^n)\cap
L^{\infty}([0, T],H^s(\Bbb R^n))$ with $f\in C^{\infty}$, $f(0)=0$
and $s\ge 0$;

Sobolev's imbedding theorem of $L^{\infty}([0, T],
H^{\f{m+6}{2(m+2)}+p_0(m)-\dl}(\Bbb R^2))\subset L^{\infty}([0,
T]\times\Bbb R^2)$ for small $\dl>0$ and $m\le 9$.

For small $T$, one can derive from (5.7) that
$$\|\Cal F(w)(t,\cdot)\|_{H^{\f{m+6}{2(m+2)}+p_0(m)-\dl}}\le\f13.\tag5.8$$

On the other hand, by (3.25) in  Lemma 3.4(taking $
s_0=p_1(m)-\f{\dl}{2}$), we have
$$
\align \|\Cal F(w)(t,\cdot)\|_{H^{\f{m+6}{2(m+2)}+ p_1(m)-\dl}}&\le
Ct^{2-\f{(m+2)p_1(m)}{2}+\f{(m+2)\dl}{4}}
\|f(t,\cdot,u_1(t,\cdot)+u_2(t,\cdot)+w(t,\cdot))-f(t,\cdot,0)\|_{H^{\f{m+6}{2(m+2)}-\f{\dl}{2}}}\\
&\le Ct^{2-\f{(m+2)p_1(m)}{2}+\f{(m+2)\dl}{4}},
\endalign$$
which yields for small $T$
$$t^{\f{(m+2)p_1(m)}{2}-2}\|\Cal F(w)(t,\cdot)\|_{H^{\f{m+6}{2(m+2)}+p_1(m)-\dl}(\Bbb
R^2)}\le\f13.\tag5.9$$

If we take $p_2=p_2(m)-\dl$ in (3.26) of Lemma 3.4, then we have for
small $T$
$$
\align \|\p_t\Cal F(w)(t,\cdot)\|_{H^{\f{3}{m+2}+p_2(m)-\dl}}&\le
CT^{1-\f{m+2}{2}(p_2(m)-\dl)}
\|f(t,\cdot,u_1(t,\cdot)+u_2(t,\cdot)+w(t,\cdot))-f(t,\cdot,0)\|_{H^{\f{m+6}{2(m+2)}-\dl}}\\
&\le CT^{\f{(m+2)\dl}{2}}\\
&\le\f13.\tag5.10
\endalign
$$

Collecting (5.8)-(5.10) yields for small $T$
$$\ds\sup_{t\in [0, T]}\||\Cal F(w)(t,\cdot)\||\le 1,\tag5.11$$
which means  $\Cal F$ maps $G$ into $G$.

Next we prove that the mapping $\Cal F$ in contractible for small
$T$.

For $w_1, w_2\in G$, due to
$f(\tau,x,u_1+u_2+w_1)-f(\tau,x,u_1+u_2+w_2)=\int_0^1f'(\tau,x,u_1+u_2+\th
w_1+(1-\th)w_2)(w_1-w_2)d\th$, then a direct computation yields as
in (5.8)-(5.10) for $t\in [0, T]$ and small $T$
$$\align &\sup_{t\in
(0,T]}|\|\Cal F(w_1)(t,\cdot)-\Cal F(w_2)(t,\cdot)\||\\
&=\sup_{t\in
(0,T]}|\|E(f(t,x,u_1+u_2+w_1(\tau,\cdot))-Ef(t,x,u_1+u_2+w_2(\tau,\cdot))\||\\
&\leq C(T^{\f{(m+2)\dl}{4}}+T^{\f{(m+2)\dl}{2}})\sup_{t\in
(0,T]}|\|w_1-w_2\||\\
&\le \f12\sup_{t\in (0,T]}|\|w_1-w_2\||.\tag5.12\endalign$$

Therefore, by the fixed point theorem and (5.11)-(5.12), we complete
the proof of Theorem 5.1.

Under the assumption $(A_1)$, we have

\vskip .2cm

{\bf Theorem 5.2.} {\it  Under the assumption $(A_1)$, there exists
a constant $T>0$ such that (1.1) has a local solution $u\in C([0,T],
H^{\f{n}{2}+\f{2}{m+2}-}(\Bbb R^n))\cap C((0, T],$
$H^{\f{n}{2}+\f{m+4}{2(m+2)}-}(\Bbb R^n))\cap C^1([0,T],
H^{\f{n}{2}-\f{m}{2(m+2)}-}(\Bbb R^n))$.}

{\bf Proof.} Since $C([0,T], H^{\f{n}{2}+\f{2}{m+2}-}(\Bbb
R^n))\subset L^{\infty}([0, T]\times\Bbb R^n)$, then Theorem 5.2 can
be shown by the same procedure as in Theorem 5.1, we omit it here.
\hfill \qed

Finally, we prove the local existence of solution to (1.1) under the
condition $(A_2)$.

\vskip .2cm

{\bf Theorem 5.3.} {\it  Under the assumption $(A_2)$, there exists
a constant $T>0$ such that (1.1) has a local solution $u\in
L^{\infty}{([0, T]\times\Bbb R^n)}\cap C([0,T],
H^{\f{m+6}{2(m+2)}-}(\Bbb R^n))\cap C((0, T],$
$H^{\f{m+3}{m+2}-}(\Bbb R^n))\cap C^1([0,T], H^{\f{1}{m+2}-}(\Bbb
R^n))$.}

\vskip .1cm

{\bf Proof.}  Let $u_1(t,x)$ satisfy
$$
\cases &\p_t^2u_1-t^m\Delta u_1=0,\qquad (t,x)\in [0,
+\infty)\times\Bbb R^n,\\
&u_1(0,x)=0,\quad \p_tu_1(0,x)=\vp(x),
\endcases\tag5.13
$$ where $\vp(x)$
satisfies the assumption $(A_2)$, then by Lemma 2.2 and Proposition
3.3 we know that for any fixed $\dl>0$ with $\dl<\ds\f{1}{2(m+2)}$
$$u_1(t,x)\in L^{\infty}([0, 1]\times\Bbb R^n)\cap C([0,1],
H^{\f{m+6}{2(m+2)}-\dl}(\Bbb R^2))\cap C((0, 1],
H^{\f{m+3}{m+2}-\dl}(\Bbb R^2))\cap C^1([0,1],
H^{\f{1}{m+2}-\dl}(\Bbb R^2)),$$ which satisfies for $t\in [0, 1]$
$$\|u_1(t,\cdot)\|_{L^{\infty}(\Bbb R^2)}+\|u_1(t,\cdot)\|_{H^{\f{m+6}{2(m+2)}-\dl}(\Bbb R^2)}+
t^{\f{m}{4}}\|u_1(t,\cdot)\|_{H^{\f{m+3}{m+2}-\dl}(\Bbb
R^2)}+\|\p_tu_1(t,\cdot)\|_{H^{\f{1}{m+2}-\dl}(\Bbb R^2)}\le C(
\dl).\tag5.14$$

Next, we establish the more regularities of $u_1(t,x)$ in the
directions $x'=(x_2, ..., x_n)$.

It is noted that for $|\al|\ge 1$
$$
\cases &\p_t^2\p_{x'}^{\al}u_1-t^m\Delta \p_{x'}^{\al}u_1=0,\qquad
(t,x)\in [0,
+\infty)\times\Bbb R^n,\\
&\p_{x'}^{\al}u_1(0,x)=0,\quad
\p_t\p_{x'}^{\al}u_1(0,x)=\p_{x'}^{\al}\vp(x).
\endcases\tag5.15
$$

This derives
$$\p_{x'}^{\al}u_1(t,x)\in L^{\infty}([0, 1]\times\Bbb R^n)\cap C([0,1],
H^{\f{m+6}{2(m+2)}-\dl}(\Bbb R^n))\cap C((0, 1],
H^{\f{m+3}{m+2}-\dl}(\Bbb R^n))\cap C^1([0,1],
H^{\f{1}{m+2}-\dl}(\Bbb R^n))$$ and satisfies for $t\in [0, 1]$
$$
\align &\|\p_{x'}^{\al}u_1(t,\cdot)\|_{L^{\infty}(\Bbb
R^2)}+\|\p_{x'}^{\al}u_1(t,\cdot)\|_{H^{\f{m+6}{2(m+2)}-\dl}(\Bbb
R^n)}+
t^{\f{m}{4}}\|\p_{x'}^{\al}u_1(t,\cdot)\|_{H^{\f{m+3}{m+2}-\dl}(\Bbb
R^n)}\\
&\qquad +\|\p_t\p_{x'}^{\al}u_1(t,\cdot)\|_{H^{\f{1}{m+2}-\dl}(\Bbb
R^n)}\le C_{\al}(\dl).\tag5.16\endalign$$

Set $v=u-u_1-u_2$, where $u_2$ is defined as in (5.3),  then we have
from (1.1)
$$
\cases &\p_t^2v-t^{m}\triangle v=f(t,x,u_1+u_2+v)-f(t,x,0),\\
&v(0,x)=\p_t v(0,x)=0.\endcases\tag5.17$$

In order to solve (1.1), it only suffices to solve (5.17). This
requires us to establish the a priori $L^{\infty}$ bound of
$\p_{x'}^{\al}v$ in (5.17) for $|\al|\le [\f{n}{2}]+1$. For this
end, motivated by Lemma 2.2 and Lemma 2.3, we should establish
$\p_{x'}^{\al+\beta}v\in L^{\infty}([0, T], H^s(\Bbb R^n))$ with
$s>\f12$ and $|\beta|\le [\f{n}{2}]+1$.

Taking $\p_{x'}^{\g}$ $(|\g|\le 2[\f{n}{2}]+2)$ on two hand sides of
(5.17) yields
$$
\cases &\p_t^2\p_{x'}^{\g}v-t^{m}\triangle \p_{x'}^{\g}v=
F_{\g}(t,x,
\p_{x'}^{\al}v)_{|\al|\le|\g|}\equiv\ds\sum_{|\beta|+l\le
|\g|}C_{\beta l}
\biggl((\p_{x'}^{\beta}f)(t,x,u_1+u_2+v)-(\p_{x'}^{\beta}f)(t,x,0)\biggr)\\
&\qquad \times \p_u^{l}f(t,x,u_1+u_2+v)\ds\Pi_{\Sb 1\le k\le l\\
\beta_1+...+\beta_l=l\endSb}\p_{x'}^{\beta_k}(u_1+u_2+v),\\
&\p_{x'}^{\g}v(0,x)=\p_t \p_{x'}^{\g}v(0,x)=0.\endcases\tag5.18$$

If $\ds\sum_{|\al|\le
[\f{n}{2}]+1}\|\p_{x'}^{\al}v\|_{L^{\infty}([0, T]\times\Bbb
R^n)}+\ds\sum_{|\g|\le
2[\f{n}{2}]+2}\|\p_{x'}^{\g}v\|_{L^{\infty}([0, T], H^s(\Bbb
R^n))}\le 2$ with $s>\f12$ and $T\le 1$, then by Lemma 2.3 and
(5.16), we have from (5.18) that for small $T$
$$
\ds\sum_{|\al|\le [\f{n}{2}]+1}\|\p_{x'}^{\al}v\|_{L^{\infty}([0,
T]\times\Bbb R^n)}\le 1.\tag5.19$$

Based on the preparations above,  we will use the fixed point
theorem to show Theorem 5.3.

For $w\in L^{\infty}([0, T]\times\Bbb R^n)\cap C([0,T],
H^{\f{m+6}{2(m+2)}+p_0(m)-\dl})\cap C((0,T],
H^{\f{m+6}{2(m+2)}+p_1(m)-\dl}) \cap C^1([0,T],$
$H^{\f{3}{m+2}+p_2(m)-\dl})$ with $\p_{x'}^{\al}w\in L^{\infty}([0,
T]\times\Bbb R^n)$ ($|\al|\le [\f{n}{2}]+1$) and $\p_{x'}^{\g}w\in
C([0,T], H^{\f{m+6}{2(m+2)}+p_0(m)-\dl})\cap C((0,T],
H^{\f{m+6}{2(m+2)}+p_1(m)-\dl})$ $\cap C^1([0,T],
H^{\f{3}{m+2}+p_2(m)-\dl})$ ($|\g|\le 2[\f{n}{2}]+2$), where the
expressions of $p_0(m), p_1(m), p_2(m)$ are given in (5.5), we
define
$$
\align |\|w(t,\cdot)\||&\equiv
\ds\sum_{|\al|=0}^{[\f{n}{2}]+1}\|\p_{x'}^{\al}w(t,x)\|_{L^{\infty}([0,T]\times\Bbb
R^n)}+
\ds\sum_{|\g|=0}^{2[\f{n}{2}]+2}\|\p_{x'}^{\g}w(t,\cdot)\|_{H^{\f{m+6}{2(m+2)}+
p_0(m)-\dl}(\Bbb R^n)}\\
&  +t^{\f{(m+2)p_1(m)}{2}-2}
\ds\sum_{|\g|=0}^{2[\f{n}{2}]+2}\|\p_{x'}^{\g}w(t,\cdot)\|_{H^{\f{m+6}{2(m+2)}+
p_1(m)-\dl}(\Bbb
R^n)}+\sum_{|\g|=0}^{2[\f{n}{2}]+2}\|\p_t\p_{x'}^{\g}w(t,\cdot)\|_{H^{\f{3}{m+2}+p_2(m)-\dl}(\Bbb
R^n)}.
\endalign
$$

A set $Q$ is defined as follows
$$
\align Q\equiv &\bigg\{w\in L^{\infty}([0, T]\times\Bbb R^n)\cap
C([0,T], H^{\f12+\f{2}{m+2}-\dl}(\Bbb R^n))\cap C((0,T],
H^{\f{m+3}{m+2}-2\dl}(\Bbb R^n))\cap C^1([0,T],
H^{\f{1}{m+2}-\dl}(\Bbb R^n)):\\
&\qquad \dsize\sup_{t\in [0,T]}|\|w(t,\cdot)\||\leq 2\bigg\}.
\endalign$$

Let us define a nonlinear mapping $\Cal F$ as follows
$$\Cal F(w)=E(f(t,x,u_1+u_2+w)-f(t,x,0)),\tag5.20$$
where the meaning of the operator $E$ is given in (5.5).

As in the proof procedure of Theorem 5.1, we now show that the
mapping $\Cal F$ is from $Q$ into itself and is contractible for
small $T$.

At first, $\Cal F(w)$ solves the following problem
$$
\cases &(\p_t^2-t^m\Delta)\Cal F(w)=f(t,x,u_1+u_2+w)-f(t,x,0),\\
&\Cal F(w)|_{t=0}=\p_t\Cal F(w)|_{t=0}=0.
\endcases
$$

By (5.19), we can derive that for small $T$
$$\ds\sum_{|\al|=0}^{[\f{n}{2}]+1}\|\p_{x'}^{\al}\Cal Fw(t,x)\|_{L^{\infty}([0,T]\times\Bbb
R^n)}\le  1.\tag5.21$$

Similar to the proof as in Theorem 5.1, one has for small $T$ and
$t\in [0, T]$
$$\align
&\ds\sum_{|\g|=0}^{2[\f{n}{2}]+2}\|\p_{x'}^{\g}w(t,\cdot)\|_{H^{\f{m+6}{2(m+2)}+
p_0(m)-\dl}(\Bbb R^n)}+t^{\f{(m+2)p_1(m)}{2}-2}
\ds\sum_{|\g|=0}^{2[\f{n}{2}]+2}\|\p_{x'}^{\g}w(t,\cdot)\|_{H^{\f{m+6}{2(m+2)}+
p_1(m)-\dl}(\Bbb
R^n)}\\
&\quad
+\sum_{|\g|=0}^{2[\f{n}{2}]+2}\|\p_t\p_{x'}^{\g}w(t,\cdot)\|_{H^{\f{3}{m+2}+p_2(m)-\dl}(\Bbb
R^n)}\le 1\tag5.22 \endalign$$ and
$$\||\Cal F(w_1)-\Cal F(w_2)\||\le\f12\||w_1-w_2\||,\tag5.24$$
where $w_1, w_2\in Q$.

Combining (5.21) with (5.22) yields
$$\||\Cal F(w)\||\le 2,$$ which means $\Cal F$ maps $Q$
into itself. Therefore, it follows from the fixed point theorem that
we complete the proof of Theorem 5.3.

\vskip 0.3 true cm \centerline{\bf $\S 6$. Proof on Theorem 1.1}
\vskip 0.3 true cm

Based on the results in $\S 2$-$\S 5$, we now start to prove Theorem
1.1. At first, under the assumptions $(A_3)$ and $(A_1)$, we
establish the following conclusions on the conormal regularities  of
the local solution $u(t,x)$ obtained in Theorem 5.2 and Theorem 5.3,
respectively.

\vskip .2cm

 {\bf Theorem 6.1.} {\it (i) For the solution $u(t,x)$ in
Theorem 5.2, we have $u(t,x)\in
I^{\infty}H^{\f{n}{2}-\f{m}{2(m+2)}-}(\G_0)$;

(ii). For the solution $u(t,x)$ in Theorem 5.3, then $u(t,x)\in
I^{\infty}H^{\f{1}{m+2}-}(\G_1^+\cup\G_1^-)$.}

\vskip .2cm

{\bf Proof.} (i) By the commutator relations in Lemma 4.2 and a
direct computation, we have from (1.1)
$$
\cases &\p_t^2U_k-t^m\Delta U_k=\ds\sum \Sb \beta_0+l_0\le k_0\\
\beta_{ij}+l_{ij}=k_{ij}\\
\sum l_0^s+\sum l_{ij}^s=l\le k\endSb C_{\beta
l}(L_0^{\beta_0}\Pi_{1\le i<j\le n}L_{ij}^{\beta_{ij}}\p_u^lf)(t,x,
u)\ds\Pi_{1\le s\le l}\bigl(L_0^{l_0^s}\Pi_{1\le i<j\le
n}L_{ij}^{l_{ij}^s}u\bigr),\\
&U_k(0,x)\in H^{\f{n}{2}+1-}(\Bbb R^n),\quad \p_tU_k(0,x)\in
H^{\f{n}{2}-}(\Bbb R^n),
\endcases\tag6.1
$$
here $U_k=\{L_0^{k_0}\Pi_{1\le i<j\le n}L_{ij}^{k_{ij}}u\}_{k_0+\sum
k_{ij}=k}$ for  $k\in\Bbb N\cup\{0\}$, and in the process of
deriving the regularities of $U_k(0,x)$ and $\p_tU_k(0,x)$ we have
used that facts of $\Pi_{1\le i, j\le n}(x_i\p_j)^{k_{ij}}\vp(x)$
$\in H^{\f{n}{2}-}(\Bbb R^n)$ and $w_1(x)w_2(x)\in
H^{\f{n}{2}-}(\Bbb R^n)$ if $w_1(x), w_2(x)\in H^{\f{n}{2}-}(\Bbb
R^n)$.

Next we use the induction method to prove
$$U_k(t,x)\in C([0, T], H^{\f{n}{2}+\f{2}{m+2}-}(\Bbb R^n))\cap C((0,
T], H^{\f{n}{2}+1-\f{m}{2(m+2)}-}(\Bbb R^n))\cap C^1([0, T],
H^{\f{n}{2}-\f{m}{2(m+2)}-}(\Bbb R^n))\tag6.2$$ which satisfies for
any small fixed $\dl>0$
$$\|U_k\|_{C([0, T], H^{\f{n}{2}+\f{2}{m+2}-\dl})}+t^{\f{m}{4}}
\|U_k(t,\cdot)\|_{H^{\f{n}{2}+1-\f{m}{2(m+2)}-\dl}}
+\|\p_tU_k\|_{C([0, T], H^{\f{n}{2}-\f{m}{2(m+2)}-\dl})}\le
C_k(\dl).\tag6.3$$

It is noted that (6.2)-(6.3) has been shown in Theorem 5.2 in the
case of $k=0$. Assume that (6.2)-(6.3) hold for the case up to $
k-1$, then one has by (6.1)
$$
\cases &\p_t^2U_k-t^m\Delta U_k-(\p_uf)(t,x,u)U_k=F_k(t,x),\\
&U_k(0,x)\in H^{\f{n}{2}+1-}(\Bbb R^n),\quad \p_tU_k(0,x)\in
H^{\f{n}{2}-}(\Bbb R^n),
\endcases\tag6.4
$$
where $F_k(t,x)\in C([0, T], H^{\f{n}{2}+\f{2}{m+2}-})$.

This, together with Proposition 3.3 and Lemma 3.4, yields
(6.2)-(6.3) in the case of $k$.

We now prove  $u(t,x)\in I^{\infty}H^{\f{n}{2}-}(\G_0)$.

It is noted that for $i=1, \cdots, n,$ by (6.3) and Remark 2.1,
$$N_2^i(h_2(t,x)\chi(\f{(m+2)|x|}{2t^{\f{m+2}{2}}})u)=N_2^i(h_2\chi)u
+\f{2}{m+2}(\f{(m+2)|x|}{2t^{\f{m+2}{2}}}-1)h_2\chi
t^{\f{m}{2}+1}\p_iu\in L^{\infty}([0, T],
H^{\f{n}{2}-\f{m}{2(m+2)}-}),$$  where the definitions of $h_2(t,x)$
and $\chi(\f{(m+2)|x|}{2t^{\f{m+2}{2}}})$ are given in Definition
4.3 . Furthermore, by (4.11) in Lemma 4.4 and (6.3), we can obtain
for any $k_i, k_0, k_{ij}\in\Bbb N\cup\{0\}$

$$(N_2^i)^{k_i}L_0^{k_0}\Pi_{1\le i<j\le n}L_{ij}^{k_{ij}}(h_2\chi u)\in
L^{\infty}([0, T], H^{\f{n}{2}-\f{m}{2(m+2)}-}).\tag6.5$$

This, together with (4.9) in Lemma 4.4, yields
$$\bar
L_i^{k_i}L_0^{k_0}\Pi_{1\le i<j\le n}L_{ij}^{k_{ij}}(h_2\chi u)\in
L^{\infty}([0, T], H^{\f{n}{2}-\f{m}{2(m+2)}-}).\tag6.6$$

In order to show $u(t,x)\in I^{\infty}H^{\f{n}{2}-}(\G_0)$, we need
to prove
$$\ds\Pi_{1\le
i\le n}\bar L_i^{k_i}L_0^{k_0}\Pi_{1\le i<j\le
n}L_{ij}^{k_{ij}}(h_2\chi u)\in L^{\infty}([0, T],
H^{\f{n}{2}-\f{m}{2(m+2)}-})$$ or equivalently
$$\ds\Pi_{1\le
i\le n}(N_2^i)^{k_i}L_0^{k_0}\Pi_{1\le i<j\le
n}L_{ij}^{k_{ij}}(h_2\chi u)\in L^{\infty}([0, T],
H^{\f{n}{2}-\f{m}{2(m+2)}-}).\tag6.7$$

For this end, by the commutator relations in Lemma 4.2 and (4.11),
it suffices to prove
$$N_2^{i_1}N_2^{i_2}\cdot\cdot\cdot N_2^{i_k}(h_2\chi u)\in
L^{\infty}([0, T], H^{\f{n}{2}-\f{m}{2(m+2)}-})\quad\text{for $1\le
i_1<i_2<\cdot\cdot\cdot<i_k\le n$ and $2\le k\le n$}\tag6.8$$ since
the proof on $N_2^{i_1}N_2^{i_2}\cdot\cdot\cdot
N_2^{i_k}L_0^{k_0}\Pi_{1\le i<j\le n}L_{ij}^{k_{ij}}(h_2\chi u)\in
L^{\infty}([0, T], H^{\f{n}{2}-\f{m}{2(m+2)}-})$ is completely
similar.

Indeed, by the expression of $N_2^i\equiv a(t,x)\p_i$ with
$a(t,x)=|x|-\f{2}{m+2}t^{\f{m+2}{2}}$ and (6.5), we have for $1\le
i\le n$
$$
\align \p_i^{2}\Big(a^{2}(t,x)h_2\chi u \Big)=& (a\p_i)^2(h_2\chi
u)+ \p_ia ( a \p_i)(h_2\chi u) +2a(\p_i^2a) h_2\chi u + 2 (\p_i a)^2 h_2\chi u \\
=& (N_2^i)^2 (h_2 \chi u) +(\p_i a) N_2^i(h_2 \chi u) +2a(\p_i^2a)
h_2\chi u + 2 (\p_i a)^2 h_2\chi u \\
& \in L^{\infty}([0, T], H^{\f{n}{2}-\f{m}{2(m+2)}-}), \tag6.9
\endalign
$$ here we use the facts of
$\ds\f{x_i}{|x|}\in H_{loc}^{\f{n}{2}-} (\Bbb R^n)$ and
$w_1(x)w_2(x)\in H^{\min\{s_1, s_2, s_1+s_2-\f{n}{2}\}-}(\Bbb R^n)$
if $w_1(x)\in H^{s_1}(\Bbb R^n)$ and $w_2(x)\in H^{s_2}(\Bbb R^n)$
with $s_1, s_2\ge 0$.

From (6.9), we have
$$\Delta \Big(a^{2}(t,x)h_2\chi u \Big)\in L^{\infty}([0, T], H^{\f{n}{2}-\f{m}{2(m+2)}-}),
\tag6.10$$
which derives by the regularity theory of second order elliptic
equation
$$\p_{ij}^2(a^2(t,x)h_2\chi u)\in L^{\infty}([0, T], H^{\f{n}{2}-\f{m}{2(m+2)}-})
\quad\text{for any $1\le i<j\le n$}\tag6.11$$
or equivalently
$$N_2^iN_2^j(h_2\chi u)\in L^{\infty}([0, T], H^{\f{n}{2}-\f{m}{2(m+2)}-})
\quad\text{for any $1\le i<j\le n$}.\tag6.12$$

Analogously, we can get for any $1\le i, k\le n$
$$\p_i^2 \Big(a^3\p_k(h_2\chi  u)\Big)\in L^{\infty}([0, T], H^{\f{n}{2}-\f{m}{2(m+2)}-})$$ and
$$\Delta (a^3\p_k(h_2\chi  u))\in L^{\infty}([0, T], H^{\f{n}{2}-\f{m}{2(m+2)}-}),$$ which derives
$$\p_{ij}^2(a^3\p_k(h_2\chi  u))\in L^{\infty}([0, T], H^{\f{n}{2}-\f{m}{2(m+2)}-})$$ and
further by (6.12),
$$N_2^iN_2^jN_2^k(h_2\chi  u)\in L^{\infty}([0, T], H^{\f{n}{2}-\f{m}{2(m+2)}-}).\tag6.13$$

By induction method, we can complete the proof on (6.8).

Consequently, we have
$$L_0^{k_0}\ds\Pi_{1\le i\le n}\bar L_i^{k_i}\ds\Pi_{1\le i<j\le
n}L_{ij}^{k_{ij}}(h_2\chi  u) \in \in L^{\infty}([0, T],
H^{\f{n}{2}-\f{m}{2(m+2)}-}).\tag6.14$$

Similarly, by (1) and (3) in Lemma 4.4 (noting that $\bar L_i$ can be expressed as a linear combination
of $\bar L_0$ and $L_{jk}$ with admissible coefficients in $\O_1$ and $\O_3$ respectively), we can arrive at
$$Z_1\cdots Z_k \Big( h_1 u \Big)\in
L^{\infty}([0, T], H^{\f{n}{2}-\f{m}{2(m+2)}-})\quad \text { for all
$Z_1,\cdots, Z_k\in \Cal{S}_1$}, $$ and
$$Z_1\cdots Z_k \Big( h_3 \chi_0(\f{(m+2)|x|}{2 t^{\f{m+2}{2}}}) u \Big) \in
L^{\infty}([0, T], H^{\f{n}{2}-\f{m}{2(m+2)}-})\quad \text { for all
$Z_1,\cdots, Z_k\in \Cal{S}_3$},$$ where the functions $h_1, h_3$
and $\chi_0$ are given in Definition 4.3.

Therefore,
$$u(t,x)\in I^{\infty}H^{\f{n}{2}-\f{m}{2(m+2)}-}(\G_0).$$

(ii) By the commutator relations in Lemma 4.3 and the equation
(1.1), we have for $k\ge 2$ and $j\ge 1$
$$
\cases &\p_t^2U_k-t^m\Delta U_k=\ds\sum \Sb \beta_0+l_0 \le k_0\\
\beta_{i}+l_{i}=k_{i}\\
\sum l_0^s+\sum l_i^s=l\le k\endSb C_{\beta l}(\bar
L_0^{\beta_0}\Pi_{2\le i\le n}R_i^{\beta_i}\p_u^lf)(t,x,
u)\ds\Pi_{1\le s\le l}\bigl(\bar L_0^{l_0^s}\Pi_{2\le i\le
n}R_i^{l_i^s}u\bigr),\\
&U_k(0,x)\in W^{1,\infty}(\Bbb R^n)\cap H^{\f{3}{2}-}(\Bbb
R^n),\quad \p_tU_k(0,x)\in L^{\infty}(\Bbb R^n)\cap
H^{\f{1}{2}-}(\Bbb R^n),
\endcases\tag6.15
$$
here $U_k=\{\bar L_0^{k_0}\Pi_{2\le i\le n}R_i^{k_i}u\}_{k_0+\sum
k_i=k}$ for  $k\in\Bbb N\cup\{0\}$, the definitions of $\bar L_0,
R_i (2\le i\le n)$ see Lemma 4.3, and
$(x_1\p_1)^{\g_1}\p_{x'}^{\g'}U_k(0,x)\in W^{1,\infty}(\Bbb R^n)\cap
H^{\f{3}{2}-}(\Bbb R^n),$ $
(x_1\p_1)^{\g_1}\p_{x'}^{\g'}\p_tU_k(0,x)\in L^{\infty}(\Bbb
R^n)\cap H^{\f{1}{2}-}(\Bbb R^n)$ for any multiple indices $(\g_0,
\g')$.

By Lemma 2.2-Lemma 2.3,  Proposition 3.3 and Lemma 3.4, we know from
(6.15) that
$$U_k(t,x)\in L^{\infty}([0, T]\times\Bbb R^n)\cap C([0,T],
H^{\f{m+6}{2(m+2)}-\dl}(\Bbb R^n))\cap C((0, T],
H^{\f{m+3}{m+2}-\dl}(\Bbb R^n))\cap C^1([0,T],
H^{\f{1}{m+2}-\dl}(\Bbb R^n))$$ and satisfies for $t\in [0, T]$
$$
\|U_k(t,\cdot)\|_{L^{\infty}(\Bbb
R^n)}+\|U_k(t,\cdot)\|_{H^{\f{m+6}{2(m+2)}-\dl}(\Bbb R^n)}+
t^{\f{m}{4}}\|U_k(t,\cdot)\|_{H^{\f{m+3}{m+2}-\dl}(\Bbb R^n)}
+\|\p_tU_k(t,\cdot)\|_{H^{\f{1}{m+2}-\dl}(\Bbb R^n)}\le
C_k(\dl).\tag6.16$$

Due to
$N_{2,\pm}\biggl(h_2(t,x_1)\chi_{\pm}(\f{(m+2)x_1}{2t^{\f{m+2}{2}}})u\biggr)
=N_{2,\pm}(h_2\chi)u+\biggl(\ds\f{(m+2)x_1}{2t^{\f{m}{2}+1}}\mp
1\biggr)h_2\chi t^{\f{m}{2}+1}\p_1u\in L^{\infty}([0, T],
H^{\f{1}{m+2}-\dl})$ by (6.16), here the functions $h_2$ and
$\chi_{\pm}$ are defined in Definition 4.5. Furthermore, applying
for the relations in (2) of Lemma 4.5 together with (6.16) yields

$$N_{2,\pm}^{k_1}\bar L_0^{k_0}\Pi_{2\le i\le n}L_i^{k_i}(h_2\chi u)\in L^{\infty}([0, T],
H^{\f{1}{m+2}-\dl}).\tag6.17$$

Analogously, by (1) and (3) in Lemma 4.5 and the same proof
procedure of (6.17), one can obtain
$$Z_1\cdots Z_ku(t,x))\in
L^{\infty}([0, T], H^{\f{1}{m+2}-\dl})\quad \text {on $W_i$
for all $Z_1,\cdots, Z_k\in \Cal{M}_i$}, \qquad i=1, 3.$$

Therefore,
$$u(t,x)\in I^{\infty}H^{\f{1}{m+2}-}(\G_1^+\cup\G_1^-).$$

We have completed the proof of Theorem 6.1. \hfill \qed

\vskip 0.2cm

Next, we start to illustrate $u\not\in C^2((0,T]\times\Bbb
R^2\setminus\G_0\cup\G_1^{\pm}\cup\G_2^{\pm})$ in (iii) of Theorem
1.1.

Especially, we assume that the problem (1.1) is a 2-D linear
degenerate equation  with Riemann discontinuous initial data as
follows
$$
\cases &\p_t^2u-t^m\Delta u=0,\qquad (t,x)\in [0,
+\infty)\times\Bbb R^2,\\
&u(0,x)=0,\quad \p_tu(0,x)=\vp_0(x),
\endcases\tag6.18
$$
where $\vp_0(x)=\cases
C_1\quad\text{for $x_1>0, x_2>0$}\\
C_2\quad\text{for $x_1<0, x_2>0$}\\
C_3\quad\text{for $x_1<0, x_2<0$}\\
C_4\quad\text{for $x_1>0, x_2<0$}\endcases$ with $C_i\not= C_j$ for
any $1\le i<j\le 4$ and $C_1+C_3-C_2-C_4\not=0$.

\vskip .2cm

{\bf Theorem 6.2.} {\it For the solution $u(t,x)$ of (6.18), then
$u(t,x)\not\in
I^kL^{\infty}_{loc}(\G_0\cup\G_1^{\pm}\cup\G_2^{\pm})$ with $k=2$.}

\vskip .2cm

{\bf Proof.} For convenience to write, we set
$\phi(t)=\ds\f{t^{k+1}}{k+1}$ and $\g=\ds\f{k}{2(k+1)}$ with
$k=\ds\f{m}{2}$. Then as in Lemma 2.4, the  solution of (6.18) can
be expressed as
$$u(t,x)=C_0t\int_0^1(1-s^2)^{-\g}V(s\phi(t),x)ds,\tag6.19$$
where $C_0>0$ is some fixed constant, and
$$V(\tau,x)=\cases
\vp_0(x)\quad \text{if $|x_1|\ge\tau$, $|x_2|\ge\tau$,}\\
\ds\f{C_1+C_2}{2}\quad \text{if $|x_1|\le\tau$, $x_2\ge\tau$};\quad
\ds\f{C_2+C_3}{2}
\quad \text{if $x_1\le-\tau$, $|x_2|\le\tau$};\\
\ds\f{C_3+C_4}{2}\quad \text{if $|x_1|\le\tau$, $x_2\le-\tau$};\quad
\ds\f{C_1+C_4}{2}\quad \text{if $x_1\ge\tau$, $|x_2|\le\tau$};\\
\ds\f{C_2+C_4}{2}\quad\text{if $x_1^2+x_2^2\ge\tau^2$, and $0\le
x_1\le\tau$, $0\le x_2\le\tau$ or $-\tau\le x_1\le 0$, $-\tau\le
x_2\le 0$;}\\
\ds\f{C_1+C_3}{2}\quad\text{if $x_1^2+x_2^2\ge\tau^2$, and $-\tau\le
x_1\le 0$, $0\le x_2\le\tau$ or $0\le x_1\le\tau$, $-\tau\le
x_2\le 0$;}\\
\ds\f{C_2+C_4}{2}+\ds\f{C_1+C_3-C_2-C_4}{2\pi}\p_{\tau}J(\tau,x)\quad
\text{if $x_1^2+x_2^2\le\tau^2$ and $x_1x_2>0$;}\\
\ds\f{C_1+C_3}{2}+\ds\f{C_2+C_4-C_1-C_3}{2\pi}\p_{\tau}J(\tau,x)\quad
\text{if $x_1^2+x_2^2\le\tau^2$ and $x_1x_2<0$;}\\
\endcases
$$
here
$$J(\tau,x)=\int_{\sqrt{x_1^2+x_2^2}}^{\tau}\ds\f{\tau}{\sqrt{\tau^2-r^2}}
\biggl(\f{\pi}{2}-arcsin(\f{x_1}{r})-arcsin(\f{x_2}{r})\biggr)rdr=\tau
J(1,\f{x}{\tau})\equiv \tau\t J(\f{x}{\tau})$$ with $$\t
J(y)=\int_{\sqrt{y_1^2+y_2^2}}^1\ds\f{1}{\sqrt{1-r^2}}
\biggl(\f{\pi}{2}-arcsin(\f{y_1}{r})-arcsin(\f{y_2}{r})\biggr)rdr.$$

\noindent Due to $$\p_1\t J(y)=arccos\ds\f{y_2}{\sqrt{1-y_1^2}}
\quad {and} \quad \p_2\t J(y)=arccos\ds\f{y_1}{\sqrt{1-y_2^2}},$$
then
$$
\cases &\p_1J=arccos\ds\f{x_2}{\sqrt{\tau^2-x_1^2}},\quad
\p_2J=arccos\ds\f{x_1}{\sqrt{\tau^2-x_2^2}},\\
&\p_{\tau}J=\ds\f{J}{\tau}
-\ds\f{x_1}{\tau}arccos\ds\f{x_2}{\sqrt{\tau^2-x_1^2}}
-\ds\f{x_2}{\tau}arccos\ds\f{x_1}{\sqrt{\tau^2-x_2^2}}.
\endcases\tag6.20
$$

It follows from (6.19), Remark 2.4 and a direct computation that for
$(t,x)\in\O_1^+\cap\{|x|<\phi(t)\}$
$$
\cases
&\p_tu=\ds\f{u}{t}+C_0\int_0^1(1-s^2)^{-\g}(k+1)\tau\p_{\tau}V(\tau,
x)|_{\tau=\phi(t)s}ds,\\
&\p_1u=C_0t\int_0^1(1-s^2)^{-\g}\p_1V(\tau,
x)|_{\tau=\phi(t)s}ds,\\
&\p_2u=C_0t\int_0^1(1-s^2)^{-\g}\p_2V(\tau,
x)|_{\tau=\phi(t)s}ds.\\
\endcases\tag6.21
$$

Thus, in domain $\O_1^+=\{(t,x): 0<t<1, |x_1-\phi(t)|<\dl\phi(t),
|x_2|<\dl\phi(t)\}$,
$$
\align &\f12
L_0u=u+(k+1)C_0t\int_0^1(1-s^2)^{-\g}(\tau\p_{\tau}V+x_1\p_1V+x_2\p_2V)|_{\tau=\phi(t)s}ds,\\
& \f12\bar
L_1u=\f{x_1}{\phi(t)}u+C_0(k+1)t\int_0^1(1-s^2)^{-\g}\biggl\{s(x_1\p_{\tau}V+\tau\p_1V)|_{\tau=\phi(t)s}
+\f{1-s^2}{s}(\tau \p_1V)|_{\tau=\phi(t)s}\biggr\}ds.\\
\endalign
$$

In addition, for $|x|<\tau$ and by (6.20), one has
$$
\cases &\p_{\tau}V(\tau,x)=-\ds\f{x_1x_2}{\sqrt{\tau^2-|x|^2}}
\biggl(\ds\f{1}{\tau^2-x_1^2}+\ds\f{1}{\tau^2-x_2^2}\biggr),\\
&\p_1V(\tau,x)=\ds\f{x_2\tau}{(\tau^2-x_1^2)\sqrt{\tau^2-|x|^2}},\\
&\p_2V(\tau,x)=\ds\f{x_1\tau}{(\tau^2-x_2^2)\sqrt{\tau^2-|x|^2}}\\
\endcases\tag6.22$$
and further for $(t,x)\in\O_1^+\cap\{(t, x): |x|<\phi(t)\}$ and $|x|
\le \tau,$
$$
\cases &\tau\p_{\tau}V+x_1\p_1V+x_2\p_2V=0,\\
&|x_1\p_{\tau}V+\tau\p_1V|=\ds\f{|x_2|\sqrt{\tau^2-|x|^2}}{\tau^2-x_2^2}\le
\f{|x_2|\sqrt{\tau^2-|x|^2}}{x_1^2}\in L^{\infty},\\
&|(x_1\p_{\tau}+\tau\p_1)^2V|=\ds\f{2x_1|x_2|\tau\sqrt{\tau^2-|x|^2}}{(\tau^2-x_2^2)^2}\in L^{\infty},\\
&\p_1(x_1\p_{\tau}V+\tau\p_1V)=-\ds\f{x_1x_2}{(\tau^2-x_2^2)\sqrt{\tau^2-|x|^2}},\\
&\p_1^2V=\ds\f{x_1x_2\tau(3\tau^2-2|x|^2-x_1^2)}{(\tau^2-x_1^2)^2(\tau^2-|x|^2)^{\f32}}.
\endcases\tag6.23$$

From (6.23), we can arrive at $L_0u=0$ and further $L_0^lu=0$ for
any $l\in\Bbb N\cup\{0\}$.

We now show $\bar L_1u\in L^{\infty}(\O_1^+\cap\{(t, x):
|x|<\phi(t)\})$ but $\bar L_1^2u\not\in L^{\infty}(\O_1^+\cap\{(t,
x): |x|<\phi(t)\})$.

By the expression of $\bar L_1u$ and (6.23), it suffices to prove
$$t\phi(t)\int_0^1(1-s^2)^{1-\g}\p_1V(s\phi(t),x)ds\in L^{\infty}(\O_1^+).\tag6.24$$

By the expression of $V(\tau, x)$, we only require to take care of
$\p_1V(s\phi(t),x)$ in the domain $\{(t,x): |x|\le \phi(t)s\}$ in
(6.24). At this time, a direct computation yields for
$(t,x)\in\O_1^+\cap\{(t,x): |x|<\phi(t)\}$
$$
|t\phi(t)\int_0^1(1-s^2)^{1-\g}\p_1V(s\phi(t),x)ds|\le
t\int_{\f{|x|}{\phi(t)}}^1\ds\f{|x_2|\phi^2(t)s}{\sqrt{\phi^2(t)s^2-|x|^2}}
\ds\f{1}{\phi^2(t)s^2-x_1^2}ds\equiv A_1(t,x).\tag6.25
$$

We can assert
$$A_1(t,x)\in L^{\infty}(\O_1^+\cap\{(t, x):
|x|<\phi(t)\}).\tag6.26$$

Indeed, if we set $a=\sqrt{\ds\f{\phi(t)-|x|}{\phi(t)+|x|}}$ and
$\xi=\sqrt{\ds\f{\phi(t)s-|x|}{\phi(t)s+|x|}}$, then $A_1(t,x)$ can
be estimated as follows
$$
\align
|A_1&(t,x)|=2t\biggl|\int_0^a\ds\f{x_2|x|(1+\xi^2)}{4x_1^2\xi^2+x_2^2(1+\xi^2)^2}d\xi\biggl|\\
&=2t|x_2||x|\biggl|\int_0^a\ds\f{1+\xi^2}{x_2^2(1+\xi^2+\f{2x_1^2+2x_1|x|}{x_2^2})(1+\xi^2
+\f{2x_1^2-2x_1|x|}{x_2^2})}d\xi\biggl|\\
&=\ds\f{t|x_2|}{2x_1}\biggl|\int_0^a\biggl(\ds\f{\f{2x_1^2+2x_1|x|}{x_2^2}}{1+\xi^2
+\f{2x_1^2+2x_1|x|}{x_2^2}}
-\ds\f{\f{2x_1^2-2x_1|x|}{x_2^2}}{1+\xi^2+\f{2x_1^2-2x_1|x|}{x_2^2}}\biggr)d\xi\biggl|\\
&=t \ \biggl|arctan\eta\Big|_{\eta=0}^{\f{ax_2}{|x|+x_1}}
+arctan\eta\Big|_{\eta=0}^{\f{ax_2}{|x|-x_1}}\biggr|\le C t.\\
\endalign
$$

Consequently, (6.26) holds true.

Next we show that
$$\bar L_1^2u\not\in L^{\infty}(\O_1^+\cap\{(t, x):
|x|<\phi(t)\}).\tag6.27$$

For $(t,x)\in \O_1^+\cap\{(t, x):
|x|<\phi(t)\}$, we write
$$
t\phi^2(t)\int_{\f{|x|}{\phi(t)}}^1
(1-s^2)^{2-\g}\p_1^2V|_{\tau=s\phi(t)}ds\equiv B_1(t,x)+B_2(t,x),
$$
where
$$
\align &B_1(t,x)=t\phi^3(t)\int_{\f{|x|}{\phi(t)}}^1s(1-s^2)^{2-\g}\ds\f{3x_1x_2}{(\phi^2(t)s^2-x_1^2)^2
\sqrt{\phi^2(t)s^2-|x|^2}}ds\\
&B_2(t,x)=t\phi^3(t)\int_{\f{|x|}{\phi(t)}}^1s(1-s^2)^{2-\g}\ds\f{x_1x_2^3}{(\phi^2(t)s^2-x_1^2)^2
(\phi^2(t)s^2-|x|^2)^{\f32}}ds\\
\endalign
$$

As in the process to treat $A_1(t,x)$, we set
$a=\sqrt{\ds\f{\phi(t)-|x|}{\phi(t)+|x|}}$ and
$\xi=\sqrt{\ds\f{\phi(t)s-|x|}{\phi(t)s+|x|}}$, then
$$
\align |B_1&(t,x)|=
x_1|x_2||x|t\phi(t)\int_0^a(1+\xi^2)\biggl(1-(\f{|x|}{\phi(t)}\f{1+\xi^2}{1-\xi^2})^2\biggr)^{2-\g}
\f{(1-\xi^2)^2}{(x_2^2+2(x_1^2+|x|^2)\xi^2+x_2^2\xi^4)^2}d\xi\\
<& +\infty,\tag6.28\\
|B_2&(t,x)|=\f{Cx_1|x_2|^3t\phi(t)}{|x|}\int_0^a(1+\xi^2)\biggl(1-(\f{|x|}{\phi(t)}\f{1+\xi^2}
{1-\xi^2})^2\biggr)^{2-\g}
\f{(1-\xi^2)^4}{(x_2^2+2(x_1^2+|x|^2)\xi^2+x_2^2\xi^4)^2}\ds\f{1}{\xi^2}d\xi\\
=& +\infty. \qquad\text{(Due to $\xi=0$ is a singularity point, the
integrand behaves like $\f{1}{\xi^2}$ near $\xi=0$)}\tag6.29
\endalign
$$

It is noted that the integrand in $B_2(t,x)>0$ is positive for $x_2>0$ or negative for $x_2<0$ respectively,
then by the definition of partial derivatives together with Fatou's lemma (i.e., for $G(y)=\int_a^b g(s,y)ds$
and $\p_yg(s,y)>0$, then $\ds{\underline\lim}_{h\to 0+}\f{G(y+h)-G(y)}{h}\ge\int_a^b \p_yg(s,y)ds$), we 
have from a direct computation that for
$(t,x)\in \O_1^+\cap\{(t, x): |x|<\phi(t)\}$ and if $u\in C^2((0,T]\times\Bbb
R^2\setminus\G_0\cup\G_1^{\pm}\cup\G_2^{\pm})$
$$
\align |\bar L_1^2u|\ge &
\biggl|Ct\phi^2(t)\bigl|\int_{\f{|x|}{\phi(t)}}^1 s(1-s^2)^{1-\g}
\p_1^2V(\tau, x)|_{\tau=s\phi(t)}ds\bigl| -\text{$L^{\infty}$
terms}\\
&\qquad -Ct\phi(t)\int_{\f{|x|}{\phi(t)}}^1s(1-s^2)^{-\g}\bigl|(\p_1(x_1\p_{\tau}
+\tau\p_1)V\bigr)|_{\tau=s\phi(t)}|ds\biggl|\\
=&\biggl|Ct\phi^2(t)\int_{\f{|x|}{\phi(t)}}^11-s^2)^{2-\g}|\p_1^2V|_{\tau=s\phi(t)}|ds-\text{$L^{\infty}$
terms}-CtA_1(t,x)\biggl|\\
\ge & |C B_1(t,x)+ C B_2(t,x)-\text{$L^{\infty}$
terms}|\qquad\text{(here we have used (6.26))}\\
=& +\infty \qquad\text{(here we have used (6.28)-(6.29))}.
\endalign
$$

Therefore, (6.27) is proved, and $u(t,x)\not\in
I^kL^{\infty}_{loc}(\G_0\cup\G_1^{\pm}\cup\G_2^{\pm})$ with $k=2$.
\hfill \qed

\vskip .2cm

Finally, we can complete the proof on Theorem 1.1.

\vskip .2cm

 {\bf Proof of Theorem 1.1.} (i) Combining Theorem 5.2
with Theorem 6.1.(i) yields its proof.

(ii)  Its proof comes from Theorem 5.3 and Theorem 6.1.(ii).

(iii)  Based on Theorem 5.1 and Theorem 6.2, the proof can be
completed.



\Refs \refstyle{C}

\ref\key 1\by A.S.Barreto\paper Interactions of conormal waves for
fully semilinear wave equations\jour J. Funct. Anal. 89, no. 2,
233-273\yr 1990\endref



\ref\key 2\by M.Beals\paper Vector fields associated with the
nonlinear interaction of progressing waves\jour Indiana Univ. Math.
J. 37, no. 3, 637-666\yr 1988
\endref

\ref\key 3\by M.Beals\paper Singularities due to cusp interactions
in nonlinear waves\jour Nonlinear hyperbolic equations and field
theory (Lake Como, 1990), 36-51, Pitman Res. Notes Math. Ser., 253,
Longman Sci. Tech., Harlow, 1992\endref

\ref\key 4\by L.Bers\paper Mathematical aspects of subsonic and
transonic gas dynamics\jour Surveys in Applied Mathematics, Vol. 3
John Wiley $\&$ Sons, Inc., New York; Chapman $\&$ Hall, Ltd.,
London 1958
\endref


\ref\key 5\by J.M.Bony\paper Propagation des singularit\'es pour les
\'equations aux d\'eriv\'ees partielles non lin\'eaires\jour Exp.
No. 22, 12 pp., \'ecole Polytech., Palaiseau, 1980\endref

\ref\key 6\by J.M.Bony\paper Interaction des singularit\'es pour les
\'equations aux d\'eriv\'ees partielles non lin\'eaires\jour
Goulaouic-Meyer-Schwartz Seminar, 1981/1982, Exp. No. II, 12 pp.,
\'ecole Polytech., Palaiseau, 1982\endref


\ref\key 7\by J.M.Bony\paper Second microlocalization and
propagation of singularities for semilinear hyperbolic
equations\jour  Hyperbolic equations and related topics
(Katata/Kyoto, 1984), 11-49, Academic Press, Boston, MA, 1986\endref


\ref\key 8\by J.Y.Chemin\paper Interaction de trois ondes dans les
\'equations semi-lin\'eaires strictement hyperboliques d'ordre
2\jour Comm. Partial Differential Equations 12, no. 11, 1203-1225
\yr 1987\endref



\ref\key 9\by Chen Shuxing\paper Multidimensional Riemann problem
for semilinear wave equations \jour Comm. Partial Differential
Equations 17, no. 5-6, 715-736 \yr 1992\endref

\ref\key 10\by M.Dreher, M. Reissig\paper Propagation of mild
singularities for semilinear weakly hyperbolic equations\jour J.
Anal. Math. 82, 233-266 \yr 2000\endref

\ref\key 11\by M.Dreher, M. Reissig\paper  Local solutions of fully
nonlinear weakly hyperbolic differential equations in Sobolev
spaces\jour Hokkaido Math. J. 27, no. 2, 337-381\yr 1998\endref

\ref\key 12\by M.Dreher, I.Witt\paper Energy estimates for weakly
hyperbolic systems of the first order\jour Commun. Contemp. Math. 7,
no. 6, 809-837\yr 2005\endref


\ref\key 13\by A.Erdelyi, W.Magnus, F.Oberhettinger,
F.G.Tricomi\jour Higher Transcendental Functions, Vol.1, New York,
Toronto and London: McGraw-Hill Book Company, 1953\endref



\ref\key 14\by Han Qing, Hong Jia-Xing, Lin Chang-Shou\paper On the
Cauchy problem of degenerate hyperbolic equations\jour Trans. Amer.
Math. Soc. 358, no. 9, 4021-4044 \yr 2006\endref

\ref\key 15\by Han Qing\paper Energy estimates for a class of
degenerate hyperbolic equations\jour Math. Ann. 347, no. 2,
339-364\yr 2010\endref


\ref\key 16\by L.H$\ddot o$rmander, The analysis of linear partial
differential operators I, Distribution theory and Fourier analysis,
Reprint of the second (1990), Classics in Mathematics.
Springer-Verlag, Berlin, 2003\endref

\ref\key 17\by  D.Lupo, C.S.Morawetz,  K.R.Payne\paper On closed
boundary value problems for equations of mixed elliptic-hyperbolic
type\jour Comm. Pure Appl. Math. 60, no. 9, 1319-1348\yr 2007\endref



\ref\key 18\by R.Melrose, N.Ritter\paper Interaction of nonlinear
progressing waves for semilinear wave equations\jour Ann. of Math.
(2) 121, no. 1, 187-213\yr 1985\endref



\ref\key 19\by G.M\'etivier\paper The Cauchy problem for semilinear
hyperbolic systems with discontinuous data\jour Duke Math. J. 53,
no. 4, 983-1011\yr 1986\endref

\ref\key 20\by G.M\'etivier, J.Rauch\paper Interaction of piecewise
smooth progressing waves for semilinear hyperbolic equations\jour
Comm. Partial Differential Equations 15, no. 8, 1079-1140 \yr
1990\endref


\ref\key 21\by  C.S.Morawetz\paper Mixed equations and transonic
flow\jour J.Hyperbolic Differ. Equ. 1, no.1, 1-26\yr 2004\endref


\ref\key 22\by  J.M.Rassias\paper Mixed type partial differential
equations with initial and boundary values in fluid mechanics\jour
Int. J. Appl. Math. Stat. 13, No. J08, 77-107\yr 2008\endref


\ref\key 23\by K.Taniguchi, Y.Tozaki\paper A hyperbolic equation
with double characteristics which has a solution with branching
singularities\jour Math.Japan 25, 279-300\yr 1980\endref


\ref\key 24\by K.Yagdjian\paper A note on the fundamental solution
for the Tricomi-type equation in the hyperbolic domain\jour J.
Differential Equations 206, no. 1, 227-252 \yr 2004\endref

\ref\key 25\by K.Yagdjian\paper Global existence for the
n-dimensional semilinear Tricomi-type equations \jour Comm. Partial
Differential Equations 31, no. 4-6, 907-944 \yr 2006\endref

\endRefs
\enddocument